\documentclass{amsart}
\usepackage{amscd,amssymb,hyperref}
\usepackage{diagrams}

\theoremstyle{plain}
\newtheorem{prop}[subsection]{Proposition}
\newtheorem{thm}[subsection]{Theorem}
\newtheorem*{theo}{Theorem}

\newtheorem{lem}[subsection]{Lemma}
\newtheorem{cor}[subsection]{Corollary}
\newtheorem{definition}[subsection]{Definition}
\theoremstyle{remark}
\newtheorem{rem}[subsection]{Remark}
\newtheorem*{ack}{Acknowledgment}

\theoremstyle{definition}
\newtheorem{exm}[subsection]{Example}

\numberwithin{equation}{section}

\renewcommand{\b}[1]{\mathbf{#1}}
\renewcommand{\atop}[2]{\genfrac{}{}{0pt}{}{#1}{#2}}

\newcommand{\A}{{\mathcal A}}
\newcommand{\Ai}{{\mathcal A}_\infty}
\newcommand{\cB}{{\mathcal B}}
\newcommand{\cG}{{\mathcal G}}
\newcommand{\LL}{{\mathcal L}}
\newcommand{\cS}{{\mathcal S}}
\newcommand{\cT}{{\mathcal T}}
\newcommand{\cV}{{\mathcal V}}
\newcommand{\cW}{{\mathcal W}}
\newcommand{\Z}{{\mathbb Z}}
\newcommand{\C}{{\mathbb C}}
\newcommand{\CP}{{\mathbb{CP}}}
\newcommand{\bl}{{\boldsymbol{\lambda}}}
\newcommand{\sfA}{{\sf A}}
\newcommand{\sfB}{{\sf B}}
\newcommand{\sfb}{{\sf b}}
\newcommand{\sfM}{{\sf M}}
\newcommand{\la}{{\lambda }}
\newcommand{\bul}{{\bullet }}
\newcommand{\p}{\partial}

\renewcommand{\c}{{\gamma }}
\renewcommand{\ll}{{\ell }}

\DeclareMathOperator{\Span}{span}
\DeclareMathOperator{\rank}{rank}
\DeclareMathOperator{\codim}{codim}

\DeclareMathOperator{\ii}{i}
\DeclareMathOperator{\im}{im}

\DeclareMathOperator{\Aut}{Aut}
\DeclareMathOperator{\Dep}{Dep}
\DeclareMathOperator{\GM}{\Omega}

\begin{document}

\title[Gauss-Manin connections for arrangements, IV]
{Gauss-Manin connections for arrangements, IV\\ Nonresonant eigenvalues}
\author[D.~Cohen]{Daniel C.~Cohen$^\dag$}
\address{Department of Mathematics, Louisiana State University,
Baton Rouge, LA 70803}
\email{\href{mailto:cohen@math.lsu.edu}{cohen@math.lsu.edu}}
\urladdr{\href{http://www.math.lsu.edu/~cohen/}
{http://www.math.lsu.edu/\~{}cohen}}
\thanks{{$^\dag$}Partially supported by National Security Agency grant
MDA904-00-1-0038}

\author[P.~Orlik]{Peter Orlik$^\ddag$}
\address{Department of Mathematics, University of Wisconsin,
Madison, WI 53706}
\email{\href{mailto:orlik@math.wisc.edu}{orlik@math.wisc.edu}}
\thanks{{$^\ddag$}Partially supported by National Security Agency
grant MDA904-02-1-0019}



\begin{abstract}
An arrangement is a finite set of hyperplanes in a finite dimensional 
complex affine space.
A complex rank one local system on the arrangement complement is
determined by a set of complex weights for the hyperplanes.
We study the Gauss-Manin connection for the moduli space of 
arrangements of fixed combinatorial type in the cohomology of  
the complement with coefficients in the local system determined by 
the  weights.  For nonresonant weights, we solve the eigenvalue problem 
for the endomorphisms arising in the $1$-form associated to the Gauss-Manin connection.
\end{abstract}

\date{\today}

\maketitle

\section{Introduction} \label{sec:intro}
Let $\A=\{H_1,\dots,H_n\}$ be an arrangement of $n$ ordered
hyperplanes in $\C^\ll$, with complement 
$\sfM=M(\A)=\C^\ll\setminus\bigcup_{j=1}^n H_j$.  
Assume that $\A$ contains $\ell$ linearly independent hyperplanes.  
A complex rank one local system on $\sfM$ is determined by a
collection of weights $\bl=(\la_1,\dots,\la_n)\in\C^n$.  Associated to
$\bl$, we have a representation $\rho:\pi_1(\sfM)\to\C^*$, given by
$\c_j\mapsto \exp(-2\pi\ii\la_j)$ for any meridian loop $\c_j$ about
the hyperplane $H_j$ of $\A$, and an associated local system
$\LL$ on $\sfM$.  
For weights which are nonresonant in the sense of Schechtman,
Terao, and Varchenko \cite{STV}, the local system cohomology vanishes in 
all but one dimension, $H^q(\sfM;\LL)=0$ for $q\neq \ell$.
Parallel translation of fibers over curves in
the moduli space of all arrangements combinatorially equivalent to $\A$ 
gives rise to a Gauss-Manin connection on the vector bundle
over this moduli space with fiber $H^\ell(\sfM;\LL)$.  
This connection arises in a variety of 
applications, including the Aomoto-Gelfand theory 
of hypergeometric integrals \cite{AK,Gel1,OT2}, 
and the representation theory of Lie algebras and quantum groups~\cite{SV,Va}.  As such, it has been studied by a number of authors, 
including Aomoto~\cite{Ao},
Schechtman and Varchenko \cite{SV,Va}, Kaneko~\cite{JK}, and Kanarek
\cite{HK}.

Denote the combinatorial type of $\A$ by $\cT$.  The moduli space 
of all arrangements of type $\cT$ 
is determined by the set of dependent collections 
of subsets of hyperplanes in the projective closure of $\A$ in $\CP^\ell$, 
see \cite{T1}.  
Let $\sfB(\cT)$ be a smooth, connected component of this
moduli space.  There is a fiber bundle $\sf{p}:\sfM(\cT) \to \sfB(\cT)$
whose fibers, $\sf{p}^{-1}(\sfb)=\sfM_{\sfb}$, are complements of
arrangements $\A_{\sfb}$ of type $\cT$.  Since $\sfB(\cT)$ is
connected, $\sfM_{\sfb}$ is diffeomorphic to $\sfM$.
The fiber bundle $\sf{p}:\sfM(\cT) \to \sfB(\cT)$ is locally trivial.
Consequently, given a local system on the fiber, there is an
associated flat 
vector bundle $\b{H}\to\sfB(\cT)$, with fiber
$H^{\ell}(\sfM_{\sfb};\LL_{\sfb})$ at
$\sfb\in\sfB(\cT)$.
For nonresonant weights, Terao \cite{T1}  
showed that the Gauss-Manin connection on this vector bundle 
has connection $1$-form  
\begin{equation} \label{eq:1form}
\nabla = \sum \Theta_{\cT'} \otimes \Omega_\bl(\cT',\cT),
\end{equation}
where $\Theta_{\cT'}$ is a logarithmic $1$-form on the closure of  
$\sfB(\cT)$ with a simple pole along the divisor 
corresponding 
to the codimension one degeneration $\cT'$ of $\cT$, and $\Omega_\bl(\cT',\cT)$ 
is an endomorphism of $H^\ell(\sfM;\LL)$.  For general position arrangements, 
this Gauss-Manin connection was found by Aomoto and Kita
\cite{AK}. Terao \cite{T1} computed this connection for
a larger class of arrangements.  In \cite{CO4}, we determined
the ``Gauss-Manin endomorphisms'' $\Omega_\bl(\cT',\cT)$ for all arrangements.  
The aim of this paper is to solve the eigenvalue problem for these 
endomorphisms.

Identify the hyperplanes of $\A$ with their indices. 
An edge of  $\A$ is a nonempty
intersection of hyperplanes in $\A$.  An edge is {\em dense} if
the subarrangement of hyperplanes containing it is irreducible: the
hyperplanes cannot be partitioned into nonempty sets so that, after a
change of coordinates, hyperplanes in different sets are in different
coordinates, see
\cite{STV}.  For an edge $X$, define $\la_X=\sum_{X \subseteq
H_j}\la_j$.  Let $\Ai=\A\cup H_{n+1}$ be the projective closure of
$\A$, the union of $\A$ and the hyperplane at infinity in $\CP^\ll$,
see \cite{OT2}.  Set $\la_{n+1} = -\sum_{j=1}^n \la_j$. Schechtman, 
Terao, and Varchenko \cite{STV}, refining work of Esnault,
Schechtman, and Viehweg \cite{ESV}, found conditions on the weights
which insure that the local system cohomology groups vanish except in
the top dimension. They proved that if
 $\sfM$ is the complement of an arrangement $\A$ in
$\C^\ll$ of combinatorial type $\cT$ with $\ell$
linearly independent hyperplanes and $\LL$ is a rank one local
system on $\sfM$ whose weights $\bl$ satisfy the condition
\[
\la_X \notin \Z_{\ge 0} \text{ for every dense edge $X$ of $\Ai$,}
\]
then $H^q(\sfM;\LL) = 0$ for $q\neq\ll$ and $\dim H^\ll(\sfM;\LL) =
|\chi(\sfM)|$, where $\chi(\sfM)$ is the Euler characteristic of
$\sfM$. These conditions depend only on the type $\cT$, so
we call weights satisfying them $\cT$-nonresonant.

Throughout this paper, we assume that $\A$ contains 
$\ell$ linearly
independent hyperplanes, hence $n \geq \ell$, and that $\bl$ is 
$\cT$-nonresonant.  We consider only codimension one degenerations 
of combinatorial types and refer to these as simply degenerations.

\begin{theo}
Let $\cT'$ be a degeneration of $\cT$, and let 
$\bl$ be a collection of generic $\cT$-nonresonant weights for the
rank one local system $\LL$.  Then the Gauss-Manin endomorphism 
$\GM_{\bl}(\cT',\cT)$ is diagonalizable.  
The spectrum of $\Omega_\bl(\cT',\cT)$ is contained in the set 
$\{0,\la_S\}$, where 
$\la_S=\sum_{j \in S}\la_j$
for some $S \subset \{1,\ldots,n+1\}$. 
\end{theo}
The set $S$ is part of a pair $(S,r)$, called 
the {\em principal dependence} of the degeneration
$\cT'$ of $\cT$, see 
Theorem \ref{thm:principal}.  
It follows from our results in Sections \ref{sec:gp}, 
\ref{sec:formal}, and \ref{sec:evalues} 
that weights $\bl$ which satisfy $\la_S\neq 0$ are sufficiently generic.  
Our results also yield an algorithm for determining the 
multiplicities of the eigenvalues, see Remark \ref{rem:algorithm}.

Let $\cG$ denote the combinatorial type of a general position 
arrangement of $n$ hyperplanes in $\C^\ll$.  The cohomology of the complement 
of an arrangement of type $\cG$ is the rank $\ell$ truncation, $A^\bul(\cG)$, 
of the exterior algebra
on $n$ generators $e_j$, $j\in [n]$,  where $[n]=\{1,\ldots,n\}$,
corresponding to the hyperplanes. 
The Orlik-Solomon algebra  $A^\bul(\A)\simeq H^\bul(M(\A);\C)$ 
is generated by one
dimensional classes $a_j$, $j\in [n]$.
It is the quotient of 
$A^\bul(\cG)$ by a homogeneous ideal, $I^\bul(\A)$,
hence it is a finite dimensional graded $\C$-algebra \cite{OT1}.
It is known that $A^\bul(\A)$ depends only on the combinatorial
type $\cT$ of $\A$ so we may write $A^\bul(\cT)$.  

Weights $\bl$ yield an element $a_{\bl}=\sum_{j=1}^n \la_j a_j$ in
$A^1(\cT)$, and multiplication by $a_{\bl}$ gives $A^\bul(\cT)$ the structure 
of
a cochain complex.  The resulting cohomology  $H^\bul(\cT)=H^\bul(A^\bul(\cT),a_\bl)$ is
a combinatorial analog of $H^\bul(M(\A);\LL)$. 
If the weights are $\cT$-nonresonant, then  $H^\bul(M(\A);\LL) \simeq 
H^\bul(A^\bul(\cT),a_\bl)$ and the only (possibly) nonzero group $H^\ell(\cT)$ 
has the $\beta${\bf nbc} basis of Falk and Terao \cite{FT}. 
This basis provides an explicit surjection 
$\tau: H^\ell(\cG) \rightarrow H^\ell(\cT)$. 
Our results in \cite{CO4} yield 
a commutative diagram of endomorphisms for each degeneration $\cT'$ of $\cT$:

\begin{equation} \label{eq:main1}
\begin{CD}
H^\ell(\cG) @>\tau>> H^\ell(\cT)\\
@VV\widetilde{\Omega}_\bl(\cT',\cT)V      
@VV\Omega_\bl(\cT',\cT)V \\
H^\ell(\cG) @>\tau>> H^\ell(\cT)
\end{CD}
\end{equation}
The endomorphism $\widetilde{\Omega}_\bl(\cT',\cT)$ of $H^\ell(\cG)$ is induced 
by an endomorphism $\omega_\bl^\bul(\cT',\cT)$ of $A^\ell(\cG)$, see \cite{GM3},  
\eqref{eq:GP endo}, and Theorem \ref{thm:induced endo}.

Here is a brief outline of the paper.
In Section \ref{sec:gp}, we recall the Aomoto complex and the 
``formal Gauss-Manin connection matrices'' of \cite{GM3} 
which are essential in our arguments.
We recall the moduli space of combinatorially equivalent arrangements
in Section \ref{sec:principal} and identify  the principal dependence $(S,r)$ of the
degeneration $\cT'$ of $\cT$.
Using the principal dependence, we construct a realizable type
$\cT(S,r)$ and an endomorphism ${\Omega}_\bl(S,r)$ of $H^\ell(\cG)$.
In Section \ref{sec:formal},   
we  determine the eigenstructure of the endomorphism
${\Omega}_\bl(S,r)$.  
In Section \ref{sec:evalues}, we show that 
$\widetilde{\Omega}_\bl(\cT',\cT)$ 
may be replaced by ${\Omega}_\bl(S,r)$ in
\eqref{eq:main1} and thereby
 determine the eigenstructure of $\Omega_\bl(\cT',\cT)$. We conclude
 with several examples to illustrate the main result.

\section{General position}\label{sec:gp}

In this section, we record a number of constructions in the Orlik-Solomon 
complex of a general position arrangement which will be used subsequently.

Let $\cG=\cG^\ell_n$ be 
the combinatorial type of 
a general position arrangement of $n$ hyperplanes 
in $\C^\ell$, where $n \ge \ell$.  The Orlik-Solomon algebra $A^\bul(\cG)$ is 
the 
rank $\ell$ truncation of an exterior algebra on $n$ generators.  Let
$T=\{i_1,\ldots, i_q\} \subset [n]$.  If order matters, we call $T$ a
$q$-tuple and write $T=(i_1,\ldots, i_q)$ and $e_T=e_{i_1} \cdots
e_{i_q}$.  The algebra $A^\bul(\cG)$ is generated (as an algebra) by 
$\{e_j\mid 
1\le j \le n\}$,  
and has (additive) basis $\{e_T \}$, where $e_T=1$ if $T=\emptyset$, and $T 
\neq \emptyset$
is an increasingly ordered tuple of 
cardinality at most~$\ell$.

Define a map $\partial:A^q(\cG) \to A^{q-1}(\cG)$ by 
$\partial(e_T)=\sum_{k=1}^q (-1)^{k-1} e_{T_k}$, 
where $T_k = (i_1,\dots,\widehat{i_k},\dots,i_q)$ if $T=(i_1,\dots,i_q)$.  
Then $\partial \circ \partial=0$, providing $A^\bul(\cG)$ with the structure 
of a 
chain complex
\begin{equation} \label{eq:gp chain}
(A^\bullet(\cG),\partial): \quad
A^0(\cG) \xleftarrow{\partial} A^1(\cG) \longleftarrow \cdots
\longleftarrow A^{\ell-1}(\cG) \xleftarrow{\partial} A^\ell(\cG)
\end{equation}
It is well known that the homology of this complex is concentrated in the top 
dimension, $H_q(A(\cG),\partial)=0$ for $q\neq \ell$.  
The dimension of the unique nontrivial homology group is $\beta(n,\ell) = \dim 
H_\ell(A^\bul(\cG),\partial)=\sum_{k=0}^\ell 
(-1)^{k}\binom{n}{k}=\binom{n-1}{\ell}$.

Weights $\bl=(\la_1,\dots,\la_n) \in \C^n$ determine an element 
$e_\bl=\sum_{j=1}^n
\la_j e_j$ in $A^1(\cG)$.  Since $A^\bul(\cG)$ is a quotient of an exterior 
algebra, we have 
$e_\bl e_\bl=0$.  Consequently, multiplication by $e_\bl$ defines a cochain 
complex
\begin{equation} \label{eq:gp cx}
(A^\bullet(\cG),e_\bl): \quad
A^0(\cG) \xrightarrow{e_\bl} A^1(\cG) \longrightarrow \cdots
\longrightarrow A^{\ell-1}(\cG) \xrightarrow{e_\bl} A^\ell(\cG)
\end{equation}
If $\bl \neq 0$, it is well known that the cohomology of this complex is 
concentrated in the top dimension, $H^q(A^\bullet(\cG),e_\bl)=0$ for $q \neq 
\ell$, 
and that $\dim H^\ell(A^\bullet(\cG),e_\bl)=\beta(n,\ell)$.

The endomorphism ${\omega}_\bl(\cT',\cT)$ of $A^\ell(\cG)$ 
which induces 
the map 
$\widetilde{\Omega}_\bl(\cT',\cT):H^\ell(\cG) \to H^\ell(\cG)$ 
of \eqref{eq:main1} is the specialization at $\bl$ of a  
``formal Gauss-Manin connection endomorphism'' 
given in \cite{GM3} and in \eqref{eq:GP endo}.  
The latter is a linear combination of endomorphisms ${\omega}_S^\bul$ of the 
Aomoto complex $(A_R^\bullet(\cG),e_{\b{y}})$ of $\cG$, a universal complex 
for the 
cohomology $H^\bul(A^\bullet(\cG),e_\bl)$.  
The Aomoto complex has terms $A_R^q(\cG) = A^q(\cG) \otimes R$, where 
$R=\C[y_1,\dots,y_n]$ is the polynomial ring, and the boundary map  is given 
by multiplication by $e_{\b{y}}=\sum_{j=1}^n y_j e_j$.

The endomorphisms ${\omega}_S^\bul$ correspond to subsets $S$ of $[n+1]$, 
the index set of the projective closure of the general position arrangement in $\CP^\ell$, 
 $\cG_\infty$.  The symmetric group $\Sigma_{n+1}$ 
on $n+1$ letters acts on $A^\bul(\cG)$ by permuting the hyperplanes 
of $\cG_\infty$, and on $R$ by permuting the variables $y_j$, where 
$y_{n+1}=-\sum_{j=1}^n y_j$.
In the basis $\{e_j \mid 1\le j \le n\}$ for the Orlik-Solomon algebra, the 
action
of $\sigma\in \Sigma_{n+1}$  is given 
by $\sigma(e_i)=e_{\sigma(i)}$ if $\sigma(n+1)=n+1$, and by 
\[
\sigma(e_i)=\begin{cases}
-e_{\sigma(n+1)}&\text{if $\sigma(i)=n+1$,}\\
e_{\sigma(i)} - e_{\sigma(n+1)}&\text{if $\sigma(i) \neq n+1$,}
\end{cases}
\]
if $\sigma(n+1)\neq n+1$.  
Denote the induced action on the Aomoto complex by 
$\phi_\sigma:A^\bullet_R(\cG) \to A^\bullet_R(\cG)$, 
\[
\phi_\sigma(e_{i_1} \cdots e_{i_p} \otimes f(y_1,\dots,y_n))=
\sigma(e_{i_1})\cdots\sigma(e_{i_p}) \otimes 
f(y_{\sigma(1)},\dots,y_{\sigma(n)}).
\]

\begin{lem} \label{lem:sym action}
For each $\sigma \in \Sigma_{n+1}$, the map $\phi_\sigma$ is a cochain 
automorphism of the Aomoto complex $(A_R^\bullet(\cG),e_{\b{y}})$.
\end{lem}
 
If $T=(i_1,\ldots,i_p)\subset [n]$ is a $p$-tuple, then
$(j,T)=(j,i_1,\ldots,i_p)$ is the $(p+1)$-tuple which adds $j$
with $1 \leq j \leq n$ to $T$ as its first entry. 
For $S=\{s_1,\dots,s_k\} \subset [n+1]$, let 
$\sigma_S$ denote 
the permutation $\big(\begin{smallmatrix}1 & 2 & \cdots & k\\
s_1 & s_2 & \cdots & s_{k}\end{smallmatrix}\big)$.  
Write $S \equiv T$ if $S$ 
and $T$ are equal sets.  

\begin{definition} \label{def:omega}
Let $T$ be a 
$p$-tuple, $S$ a $q+1$ element subset of $[n+1]$, and $j\in[n]$. 
If $S=S_0=[q+1]$, define the endomorphism
${\omega}^\bullet_{S_0}:(A^{\bul}_R(\cG),e_{\b{y}}) \rightarrow
(A^{\bul}_R(\cG),e_{\b{y}})$ by
\[
{\omega}_{S_0}^p(e_T)=
\begin{cases}
y_j \p e_{(j,T)} & \text{if $p=q$ and $S_0\equiv (j,T)$,}\\
e_{\b{y}} \p e_T & \text{if $p=q+1$ and $S_0\equiv T$,}\\
0 & \text{otherwise.}
\end{cases}
\]
If $S\neq S_0$, 
define ${\omega}^\bullet_S = 
\phi_{\sigma_S}^{} \circ {\omega}_{S_0}^{\bullet} 
\circ 
\phi_{\sigma_S}^{-1}$.
\end{definition}
One can check that 
this agrees with the case by case definition  
in \cite[Def.~4.1]{GM3}.

\begin{prop}[{\cite[Prop.~4.2]{GM3}}] \label{prop:chain map}
For every subset $S$ of $[n+1]$, the map ${\omega}^\bullet_S$ is a cochain
homomorphism of the Aomoto complex
$(A^{\bul}_R(\cG),e_{\b{y}}))$.
\end{prop}

For $S_0=[s]$, $1\le q \le \ell$, and $1\le r \le \min(q,s-1)$, 
consider the sets 
$\cV_{S_0}^{q,r}$ and $\cW_{S_0}^{q,r}$ 
of elements in $A^q_R(\cG)$ given by 
\[ 
\begin{aligned}
\cV^{q,r}_{S_0} &= \{e_Je_K \mid |J| \le r-1\}
\bigcup \{\eta_{S_0} e_J e_K \mid|J| = r-1\}  \quad \text{and}  \\
\cW^{q,r}_{S_0} &= \{e_{S_0} e_K \ (\text{if}\ q \ge s)\} 
\bigcup \{(\partial{e_J}) e_K \mid |J| \ge r+1\} \bigcup
\{\eta_{S_0} e_J e_K \mid|J| \ge r\}, 
\end{aligned}
\]
where $J \subset S_0$, $K \subset [n] \setminus S_0$, 
and $\eta_S=\sum_{i\in S}y_i e_i$.  
Let $\cB^{q,r}_{S_0} = \cV^{q,r}_{S_0} \bigcup
\cW^{q,r}_{S_0}$.  
If 
$S \subset [n+1]$ and $S \neq S_0$, define 
$\cB_{S}^{q,r} = \{\phi_{\sigma_S}(v) \mid v \in \cB_{S_0}^{q,r}\}$.  Define $\cV_S^{q,r}$ and 
$\cW_S^{q,r}$ analogously.  
Given weights $\bl=(\la_1,\dots,\la_n)$, let  
$\cB^{q,r}_S(\bl) = \{ v|_{y_i \mapsto \la_i}
\mid v \in \cB_{S}^{q,r}\}$ 
denote 
the specialization of $\cB_S^{q,r}$ at $\bl$, a 
sets of vectors in $A^q(\cG)$.  
Define $\cV_S^{q,r}(\bl)$ and 
$\cW_S^{q,r}(\bl)$ analogously.  
We will abuse notation and write 
$\eta_S=\sum_{i\in S} \la_i e_i$ when 
working in the Orlik-Solomon algebra.  
Note that $\partial\eta_S = \la_S=\sum_{i\in S} \la_i$.

\begin{lem} \label{lem:span}  
If $\la_S \neq 0$, 
the set of vectors $\cB^{q,r}_S(\bl)$ spans
the vector space 
$A^q(\cG)$.
\end{lem}
\begin{proof}
It suffices to consider the case $S=S_0$.  

First, we show that the set $\{\partial{e}_J \mid |J|=r+1\} \bigcup 
\{\eta_S e_J \mid |J|=r-1\}$ spans $A^r(\cG^s_s)$, where 
$\cG^s_s$ is a general position arrangement of $s$ hyperplanes (indexed by 
$S$) in $\C^s$.  
For this arrangement, both the chain complex $(A^\bullet(\cG^s_s),\partial)$ 
of \eqref{eq:gp chain} and the cochain complex 
$(A^\bullet(\cG^s_s),e_\bl)=(A^\bullet(\cG^s_s),\eta_S)$ of 
\eqref{eq:gp cx} are acyclic, and
\[
\begin{aligned}
\dim \im[\partial:A^{r+1}(\cG^s_s) \to A^r(\cG^s_s)] &= 
\beta(s,r)=\binom{s-1}{r}, \\
\dim \im[\eta_S:A^{r-1}(\cG^s_s) \to A^r(\cG^s_s)]  &= 
\binom{s}{r}-\beta(s,r)=\binom{s-1}{r-1}.
\end{aligned}
\]
Note that $\dim A^r(\cG^s_s)=\binom{s}{r}=\binom{s-1}{r}+\binom{s-1}{r-1}$.

Suppose $x \in \Span\{\partial{e}_J \mid |J|=r+1\}\cap \,
\Span\{\eta_S e_J \mid |J|=r-1\}$.  Then $\partial{x}=0$, and 
$x=\eta_S y$ for some $y \in A^{r-1}(\cG_s^s)$.  So 
$\partial{x}=\partial(\eta_S y) = \la_S y-\eta_S \partial y = 0$.  
Since $\la_S\neq 0$, we can write $y=c \eta_S \partial{y}$,  
where $c=1/\la_S$.  
But this implies that $x=\eta_S y = 0$.  
Consequently, $\{\partial{e}_J \mid |J|=r+1\} \bigcup 
\{\eta_S e_J \mid |J|=r-1\}$ spans $A^r(\cG^s_s)$.

Using this, a straightforward exercise shows that the set of vectors 
$\cB^{q,r}_S(\bl)$ spans the vector space $A^q(\cG)=A^q(\cG^\ell_n)$.
\end{proof}

\section{Principal dependence}
\label{sec:principal}

Let $\cT$ be the combinatorial type of the arrangement $\A$ of $n$
hyperplanes in $\C^\ell$ with $n \ge \ell \geq 1$.  We consider
the family of all arrangements of type $\cT$.  Recall  that $\A$
is ordered by the subscripts of its hyperplanes and we assume that
$\A$, and hence every arrangement of type $\cT$, contains $\ell$
linearly independent hyperplanes.

Choose coordinates $\b{u}=(u_1,\dots,u_\ll)$ on $\C^\ll$.  The
hyperplanes of an arrangement of type $\cT$ are defined by linear
polynomials $\alpha_{i} = b_{i,0} + \sum_{j=1}^{\ell} b_{i,j} u_{j} \,\,(i
= 1,\dots, n)$.  We embed the arrangement in projective space and add
the hyperplane at infinity as last in the ordering, $H_{n+1}$.  The
moduli space of all arrangements of type $\cT$ may be viewed as the
set of matrices
\begin{equation}\label{eq:point}
\sfb=
\begin{pmatrix}
b_{1,0} & b_{1,1} & \cdots & b_{1,\ell}\\
b_{2,0} & b_{2,1} & \cdots & b_{2,\ell}\\
\vdots  & \vdots  & \ddots & \vdots \\
b_{n,0} & b_{n,1} & \cdots & b_{n,\ell}\\
1       & 0       & \cdots & 0
\end{pmatrix}
\end{equation}
whose rows are elements of $\CP^\ell$, and whose
$(\ell+1)\times(\ell+1)$ minors satisfy certain dependency
conditions, see \cite[Prop.~ 9.2.2]{OT2}.

Given $S \subset [n+1]$, let $N_S(\cT)=N_S(\sfb)$ denote the
submatrix of (\ref{eq:point}) with rows specified by $S$. Let
$\rank N_S(\cT)$ be the size of the largest minor with nonzero
determinant. Define the multiplicity of $S$ in $\cT$ by
\begin{equation} \label{eq:mult}
m_S(\cT)=|S|-\rank N_S(\cT).
\end{equation}
Call $S$ \emph{dependent} (in type $\cT$) if 
 $m_S(\cT)>0$.  
For such $S$, the linear polynomials $\{\alpha_j \mid j \in S\}$ are 
dependent.
For $q \le n+1$, let $\Dep(\cT)_q$ denote the 
dependent sets of cardinality $q$, and let $\Dep(\cT)=
\bigcup_{q} \Dep(\cT)_q$.  
If $\cT'$ is a combinatorial type 
for which $\Dep(\cT) \subset \Dep(\cT')$, let
$\Dep(\cT',\cT) = \Dep(\cT') \setminus \Dep(\cT)$.
Terao \cite{T1} showed that the combinatorial type $\cT$
is determined by $\Dep(\cT)_{\ell +1}$, but dependent sets
of both smaller and larger cardinality arise in our considerations,
see Example \ref{ex:codim1}.

Let
\[
\Dep(\cT)_q^*=
\{S\in  \Dep(\cT)_q \mid \bigcap_{j \in S}H_j \neq \emptyset\}
\] 
and
let $\Dep(\cT)^*= \bigcup_{q} \Dep(\cT)_q^*$.  If $S \in  \Dep(\cT)^*$,  then
$\codim(\bigcap_{j \in S}H_j) < |S|$. 
If $\cT'$ is a combinatorial type 
for which $\Dep(\cT)^* \subset \Dep(\cT')^*$, let
$\Dep(\cT',\cT)^* = \Dep(\cT')^* \setminus \Dep(\cT)^*$.
If $|S|\geq \ell +2$, then $S \in\Dep(\cT)$ but  $S \in\Dep(\cT)^*$
if and only if every subset of $S$ of cardinality $\ell+1$ is dependent.
It is convenient to work with these smaller collections of dependent sets.

Define endomorphisms of $A^{\bul}_{R}(\cG)$ by
\begin{equation} \label{eq:GP endo}
{\omega}^\bullet(\cT)=\sum_{S\in
\Dep(\cT)}m_S(\cT)\cdot {\omega}^\bullet_S\text{~~~~and~~~~}
{\omega}^\bullet(\cT',\cT)=\sum_{S\in
\Dep(\cT',\cT)}m_S(\cT')\cdot {\omega}^\bullet_S.
\end{equation}
These are cochain homomorphisms of the Aomoto complex by Proposition 
\ref{prop:chain map}.
Since $\Dep(\cT)_q^*= \Dep(\cT)_q$ for $q\leq \ell+1$, we have
\begin{equation} \label{eq:GP endo*}
{\omega}^\bullet(\cT)=\sum_{S\in
\Dep(\cT)^*}m_S(\cT)\cdot {\omega}^\bullet_S\text{~~~~and~~~~}
{\omega}^\bullet(\cT',\cT)=\sum_{S\in
\Dep(\cT',\cT)^*}m_S(\cT')\cdot {\omega}^\bullet_S.
\end{equation}

\begin{thm}[\cite{GM3}] \label{thm:induced endo}
The endomorphism $\widetilde{\Omega}_\bl(\cT',\cT)$ is induced by the specialization 
${\omega}_\bl(\cT',\cT):={\omega}^\ell(\cT',\cT)|_{y_j\mapsto\la_j}$ of the 
endomorphism ${\omega}^\ell(\cT',\cT)$.
\end{thm}

Denote the cardinality of $S$ by $s=|S|$.  For 
$1 \le r \le \min(\ell,s-1)$, 
consider the combinatorial type $\cT(S,r)$ defined by
\[
T \in \Dep(\cT(S,r))^* \iff |T \cap S| \geq r+1.
\]
This type is realized by a pencil of hyperplanes indexed by $S$ with a common
subspace of codimension $r$, together with $n-s$ hyperplanes in general 
position.
Note that for $r=1$ the hyperplanes in $S$ coincide, so $\cT(S,r)$   
is a multi-arrangement.

\begin{thm}
\label{thm:principal}
Let $\cT'$ be a degeneration of a realizable combinatorial type $\cT$.
For each set $S_i \in \Dep(\cT',\cT)^*$, let $r_i$ be minimal 
so that $\Dep(\cT(S_i,r_i))^* \subset  \Dep(\cT')^*$. 
Given the collection 
$\{(S_i, r_i)\}$ there is a unique pair $(S,r)$ with $r=\min \{r_i\}$, 
 $\Dep(\cT(S,r))^* \subset  \Dep(\cT')^*$, and 
 for every  pair $(S_i,r_i)$ where $r_i=r$,  
$S_i \subset S$.
\end{thm}
\begin{proof}
Terao \cite{T1} classified the three codimension one degeneration
types in the moduli space of an arrangement whose only dependent set
is the minimally dependent set $T$ of size $q+1$.  
 \begin{enumerate}
\item[I:]  $|S\cap T|\leq q-1$ for all $S \in \Dep(\cT',\cT)^*$;
\item[II:] $\{(m,T_k) \mid m \not\in T\}$ for each fixed $k$, $1\le k \le 
|T|$;
\item[III:] $ \{(m,T_k) \mid  1 \le k \le |T| \}$ for each fixed $m \not\in T$.
\end{enumerate}
If $q=1$, then Type II does not appear.  
Recall that $T_k = (i_1,\dots,\widehat{i_k},\dots,i_{q+1})$ if 
$T=(i_1,\dots,i_{q+1})$, and note that $m\in [n+1]$ in cases II and III above.

It follows from our analysis of the corresponding types in general \cite{GM3}
that if a Type II degeneration is present, then the value of $r$ decreases and
there is a unique set of maximal cardinality with minimal $r$. 
In the other types, $r$ remains constant, but a unique dependent set of
$\cT$ increases in $\cT'$.
\end{proof}

\begin{definition}\label{def:principal}
Let $\cT'$ be a degeneration of $\cT$.
We call the  pair $(S,r)$ 
which satisfies the conditions of Theorem
\ref{thm:principal} the {\em principal  dependence}
of the degeneration.
\end{definition}

\begin{exm} \label{ex:codim1}
Let $\cT$ be the combinatorial type of the arrangement $\A$ of $4$
lines in $\C^2$ depicted in Figure \ref{fig:123}.  Here 
$\Dep(\cT)^*=\{123\}$.

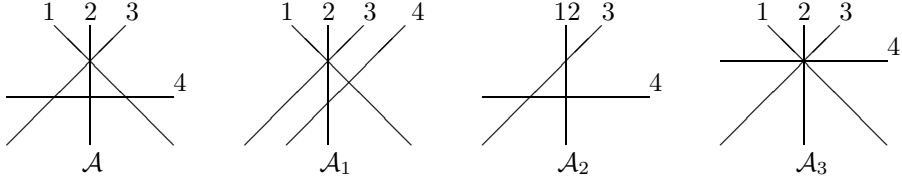
\begin{figure}[h]
\setlength{\unitlength}{.45pt}
\begin{picture}(300,130)(-200,-110)
\put(-350,-100){\line(0,1){100}} \put(-420,-100){\line(1,1){100}}
\put(-280,-100){\line(-1,1){100}}\put(-420,-60){\line(1,0){140}}
\put(-320,5){3}\put(-280,-55){4}\put(-390,5){1}\put(-355,5){2}
\put(-357,-122){$\mathcal A$}

\put(-150,-100){\line(0,1){100}}\put(-220,-100){\line(1,1){100}}
\put(-80,-100){\line(-1,1){100}}\put(-185,-100){\line(1,1){100}}
\put(-120,5){3}\put(-80,5){4}\put(-190,5){1}\put(-155,5){2}
\put(-157,-122){${\mathcal A}_{1}$}

\put(250,-100){\line(0,1){100}}\put(180,-100){\line(1,1){100}}
\put(320,-100){\line(-1,1){100}}\put(180,-30){\line(1,0){140}}
\put(280,5){3}\put(320,-25){4}\put(210,5){1}\put(245,5){2}
\put(243,-122){${\mathcal A}_{3}$}

\put(50,-100){\line(0,1){100}}\put(-20,-100){\line(1,1){100}}
\put(-20,-60){\line(1,0){140}}
\put(80,5){3}\put(120,-55){4}\put(40,5){12}
\put(43,-122){${\mathcal A}_{2}$}
\end{picture}
\caption{A line arrangement and three degenerations}
\label{fig:123}
\end{figure}

The combinatorial types $\cT_i$ of the (multi)-arrangements $\A_i$
shown in Figure~\ref{fig:123} are degenerations of
$\cT$.  For these degenerations, the  collections $\{(S_i,r_i)\}$
 and corresponding principal dependencies $(S,r)$
are given in the table below. 
\begin{center}
\begin{tabular}{l||l|l}
 &$\{(S_i,r_i)\}$&$(S,r)$\\
 \hline
$\cT_1$&$(345,2)$&$(345,2)$ \\
$\cT_2$&$(12,1), (124,2), (125,2)$ &$(12,1)$\\
$\cT_3$&$(124,2),(134,2), (234,2),(1234,2)$&$(1234,2)$ 
\end{tabular}
\end{center}
\end{exm}

For the combinatorial type $\cT(S,r)$, write 
${\omega}^\bullet(S,r)={\omega}^\bullet(\cT(S,r))$, 
see \eqref{eq:GP endo*}.
In Theorem \ref{thm:decomp} below, we show  
that the Gauss-Manin endomorphism 
$\GM_\la(\cT',\cT)$ of \eqref{eq:1form} 
is induced by the specialization 
of ${\omega}^\ell(S,r)$ at $\cT$-nonresonant weights $\bl$, 
${\omega}^\ell_\bl(S,r)$.  
First, we solve the eigenvalue problem for 
the latter endomorphism. 

\section{Diagonalization}
\label{sec:formal}

The purpose of this section is to solve the eigenvalue problem for 
${\omega}^q_\bl(S,r)$, the endomorphism of the Orlik-Solomon algebra  
obtained by specializing ${\omega}^q(S,r)$ at generic weights 
$\bl=(\la_1,\dots,\la_n)$. This allows calculation
of the eigenstructure of  the induced endomorphism in cohomology,
${\Omega}_\bl(S,r)$, 
which is related to the Gauss-Manin endomorphism in 
Theorem~\ref{thm:decomp}.
First, we establish several technical results 
concerning the endomorphism $\omega^q(S,r)$ of the Aomoto complex itself.
Recall that these endomorphisms are given explicitly by
\begin{equation*} \label{eq:(S,r)}
\omega^\bullet(S,r) = \sum_{K\in\Dep(\cT(S,r))^*} m_K(S,r)\cdot 
\omega^\bullet_K,
\end{equation*}
where $m_K(S,r)$ is the multiplicity of $K$ in type $\cT(S,r)$, see \eqref{eq:mult}, 
and $\omega^\bullet_K$ is given in Definition \ref{def:omega}.   
It follows from Proposition \ref{prop:chain map} that
$\omega^\bullet(S,r)$ is a chain map.  Note that $\omega^q(S,r)=0$ for $q < r$.

Given $(S,r)$, define
\begin{equation*}
\label{eq:Psi}
\Psi^q_{S,r} = \sum_{\atop{T \subset S}{|T|=r+1}} \omega^q(T,r).
\end{equation*}
Note that $\Psi^r_{S,r} = \omega^r(S,r)$.  
For $q\ge r$, the endomorphisms $\omega^q(S,r)$ satisfy the following 
recursion.

\begin{lem} \label{lem:recursion}
For $q \ge r$, we have
\[
\Psi^q_{S,r} = \sum_{k=0}^{s-r-1} \binom{r+k-1}{k} \omega^q(S,r+k).
\]
\end{lem}
\begin{proof}
If $T\subset [n]$ satisfies $|T|=r+1$, then $\Dep(\cT(T,r))^*=\{K \mid K 
\supseteq T\}$, 
and it is readily checked that $m_K(T,r)=1$ for each such $K$.  Hence, 
$\omega^q(T,r)=\sum_{K \supseteq T} \omega_K^q$, and we have
\[
\Psi^q_{S,r} = \sum_{\atop{T \subset S}{|T|=r+1}} \sum_{K \supseteq T} 
\omega_K^q
=\sum_{|K \cap S| \ge r+1} \omega_K^q.
\]
If $|K\cap S|=r+p$, then $\omega_K^q$ occurs $\binom{r+p}{r+1}$ times in this 
sum, so
\[
\Psi^q_{S,r} =\sum_{|K \cap S| \ge r+1} \omega_K^q
=\sum_{p\ge 1} \sum_{|K\cap S|=r+p} \binom{r+p}{r+1} \omega_K^q.
\]

If $K \in \Dep(\cT(S,j))^*$, then $|K\cap S| \ge j+1$, and 
$m_K(S,j)=|K\cap S|-j$.  It follows that $\omega^q(S,j)=\sum_{|K\cap S|\ge 
j+1}
(|K\cap S|-j)\omega_K^q$.  Hence, 
\[
\sum_{k=0}^{s-r-1} \binom{r+k-1}{k} \omega^q(S,r+k)
=\sum_{k=0}^{s-r-1} \sum_{i\ge k+1} \sum_{|K\cap S|=r+i}\binom{r+k-1}{k} (i-k) 
\omega^q_K
\] 
Rewriting this last sum, we obtain
\[
\sum_{k=0}^{s-r-1} \binom{r+k-1}{k} \omega^q(S,r+k)=
\sum_{p\ge 1} \sum_{|K\cap S|=r+p} \sum_{j=0}^{p-1} \binom{r+j-1}{j}(p-j) 
\omega^q_K.
\]
A straightforward inductive argument shows that $\sum_{j=0}^{p-1} 
\binom{r+j-1}{j}(p-j)=\binom{r+p}{r+1}$, which completes the proof.
\end{proof}

Given $S$, recall that 
$\eta_S=\sum_{i\in S} y_i e_i$ and $y_S=\sum_{i\in S}y_i
=\partial \eta_S$.

\begin{lem} \label{lem:Psi values}
Let $J \subset S$ and $L \subset[n]\setminus S$.  Then 
\[
\Psi^q_{S,r}(e_J e_K) = 
\begin{cases}
0 & \text{if $|J| \le r-1$,}\\
\binom{r+p}{r} y_S e_J e_K - \binom{r+p-1}{r-1} \eta_S (\partial{e}_J) e_K
& \text{if $|J|=r+p$, where $p\ge 0$.}
\end{cases}
\]
\end{lem}
\begin{proof}
Given $(J,L)$, it follows from Definition \ref{def:omega} that $\omega^q_K(e_J e_L) \neq 0$ 
only for the 
following $K$:
\begin{equation} \label{eq:relevant}
\begin{matrix}
(J,L),\hfill & (J_k,L,n+1), & (J,L_k,n+1), \hfill \\ (i,J,L), & (J,L,n+1), 
\hfill & (i,J_k,L,n+1), & (i,J,L_k,n+1),
\end{matrix}
\end{equation}
where $i \notin (J,L)$.

If $|J|\le r-1$, then $|K\cap S| \le r$ for each of the above $K$,  
so $T \not\subset K$ for all $T\subset S$ with $|T|=r+1$.  
It follows that $\omega^q(T,r)(e_J e_L)=0$ for each such $T$.  
Consequently, $\Psi^q_{S,r}(e_J e_L)=0$.

Let $T\subset S$ be a subset of cardinality 
$r+1$, and note that $\Psi^q_{T,r}= \sum_{K\supset T} \omega^q_K$, so 
$\Psi^q_{S,r}=\sum_{T\subset S} \Psi^q_{T,r}$, where the sum is over all 
$T\subset S$ with $|T|=r+1$.  Given such a $T$, if $|T\cap J| \le r-1$, then 
none of the sets $K$ recorded in \eqref{eq:relevant} contains $T$.  It follows 
that $\Psi^q_{T,r}(e_J e_L)=0$ if $|T\cap J| \le r-1$.

Suppose $|J|=r$.  If $|J\cap T|=r$, then $T\equiv (i,J)$ for some $i\in S 
\setminus J$, and
\[
\begin{aligned}
\Psi^q_{T,r}(e_J e_L)&= \omega_{(i,J,L)}(e_J e_L) + \sum_{k=1}^{q-r} 
\omega_{(i,J,L_k,n+1)}
(e_J e_L)\\
&=y_i \partial(e_i e_J e_L)+\sum_{k=1}^{q-r} (-1)^{r+k} y_i e_i e_J e_{L_k} \\
&=y_i e_J e_L - y_i e_i \partial(e_J e_L) + (-1)^r y_i e_i e_J \partial e_L 
=y_i e_J e_L - y_i e_i (\partial e_J) e_L.
\end{aligned}
\]
Therefore, using the identity $y_J e_J = \eta_J \partial e_J$, we have
\[
\begin{aligned}
\Psi^q_{S,r}(e_J e_L)& = \sum_{T\subset S} \Psi^q_{T,r}(e_J e_L)
=\sum_{i \in S \setminus J} (y_i e_J e_L - y_i e_i (\partial e_J) e_L) \\
&=(y_S-y_J)e_J e_L - (\eta_S-\eta_J)(\partial e_J) e_L = y_S e_J e_L - \eta_S 
(\partial e_J) e_L.
\end{aligned}
\]

Now, assume that $|J|=r+p$ for some $p\ge 1$.  
As above, we have $\Psi^q_{T,r}(e_J e_L)=0$ if $|T \cap J| \neq r,r+1$.
If $|T\cap J|=r+1$, then $T \subseteq J$ and all of the sets $K$ of 
\eqref{eq:relevant} contain $T$.    
In this instance, $\Psi^q_{T,r}(e_J e_L) = \psi(e_J e_L)$, where
\[
\psi = \omega_{(J,L)} + \omega_{(J,L,n+1)} + 
\sum_{k=1}^{q}\Bigl( \omega_{((J,L)_k,n+1)} + 
\sum_{i \notin (J,L)} \omega_{(i,(J,L)_k,n+1)}\Bigr) \\
\]
Writing $J \equiv (T,J')$, a calculation reveals that 
$\Psi^q_{T,r}(e_J e_L) = \psi(e_J e_L) = y_T e_J e_L$.

If $|T \cap J|=r$, then $T \setminus T\cap J = \{t\}$ for some $t \in 
S\setminus J$.  For such $T$, of the sets $K$ from \eqref{eq:relevant}, 
only $(t,J,L)$, $(t,J_k,L,n+1)$ for $j_k \notin T$, and $(t,J,L_k,n+1)$ 
contain $T$.  This observation, and a calculation, yields
\[
\begin{aligned}
\Psi^q_{T,r}(e_J e_L) &= \Bigl(\omega_{(t,J,L)}
+\sum_{j_k \notin T} \omega_{(t,J_k,L,n+1)} +
\sum_{k=1}^{q-r-p} \omega_{(t,J,L_k,n+1)} \Bigr)(e_J e_L)\\
&= \Bigl( \omega_{(t,J,L)} + \sum_{k=1}^{r+p} \omega_{(t,(J,L)_k,n+1)}
-\sum_{j_k \in T} \omega_{(t,J_k,L,n+1)}\bigr)(e_j e_L)\\
&= y_t e_J e_L - \sum_{j_k \in T} (-1)^{k-1} y_t e_t e_{J_k} e_L.
\end{aligned}
\]
Summing over all $T \subset S$ with $|T \cap J|=r$, we obtain
\[
\begin{aligned}
\sum_{|T\cap J|=r} &\Psi^q_{T,r}(e_J e_L) = 
\sum_{t \in S \setminus J} \sum_{\atop{A\subset[r+p]}{|A|=r}}
\bigl( y_t e_J e_L - \sum_{i=1}^r (-1)^{a_i-1} y_t e_t e_{J_{a_i}} e_L \bigr)\\
&= \binom{r+p}{r}(y_S-y_J)e_J e_L
-\sum_{t\in S\setminus J} \sum_{k=1}^{r+p} (-1)^{k-1} \binom{r+p-1}{r-1}
y_t e_t e_{J_k} e_L\\
&= \binom{r+p}{r}(y_S-y_J)e_J e_L -\sum_{t \in S\setminus J} \binom{r+p-1}{r-1}
y_t e_t(\partial e_J)e_L\\
&= \binom{r+p}{r}y_S e_J e_L - \binom{r+p-1}{r}y_J e_J e_L
-\binom{r+p-1}{r-1}\eta_S (\partial e_J) e_L.
\end{aligned}
\]
Recall that $\Psi^q_{T,r}(e_J e_L)=y_T e_J e_L$ for $T \subset J$.  
Summing over all $T\subset J$, we obtain 
$\sum_{T\subset J}\Psi^q_{T,r}(e_J e_L) = \binom{r+p-1}{r} y_J e_J e_L$. 
Therefore,
\[
\begin{aligned}
\Psi^q_{S,r}(e_J e_L)& = 
\Bigl(\sum_{T\subset S} \Psi^q_{T,r}\Bigr)(e_J e_L) =
\Bigl(\sum_{|T|=r} \Psi^q_{T,r} + \sum_{|T|=r+1} \Psi^q_{T,r}\Bigr)
(e_J e_L) \\
&= \binom{r+p}{r}y_S e_J e_L 
-\binom{r+p-1}{r-1}\eta_S (\partial e_J) e_L
\end{aligned}
\]
if $|J|=r+p$.
\end{proof}

Let $\bl=(\la_1,\dots,\la_n)$ be a collection of weights, 
and consider the endomorphism $\omega^q_\bl(S,r):A^q(\cG)
\to A^q(\cG)$ of the Orlik-Solomon 
algebra obtained by specializing $\omega^q(S,r)$ at $\bl$.  
Given $S$, 
we abuse notation and write $\eta_S=\sum_{i\in S} \la_i e_i$.  
Recall the spanning set $\cB^{q,r}_S(\bl)= \cV^{q,r}_S(\bl)
\bigcup \cW^{q,r}_S(\bl)$ of $A^q(\cG)$ from Lemma \ref{lem:span}.  

\begin{thm} \label{thm:diag}
Let $\bl$ be a collection of weights satisfying $\la_S \neq 0$.  
Then the 
specialization, ${\omega}^q_\bl(S,r)$, of ${\omega}^q(S,r)$ at 
$\bl$ is diagonalizable, with eigenvalues $0$ and $\la_S$.  
\begin{enumerate}
\item The $0$-eigenspace is spanned by the 
set of vectors $\cV^{q,r}_S(\bl)$
and has dimension 
\[
\sum_{p=0}^r \binom{s}{p}\binom{n-s}{q-p}-\binom{s-1}{r}\binom{n-s}{q-r}.
\]
\item The $\la_S$-eigenspace is spanned by 
the set of vectors $\cW^{q,r}_S(\bl)$
and has dimension
\[
\sum_{p=r+1}^{\min(q,s)} 
\binom{s}{p}\binom{n-s}{q-p}+\binom{s-1}{r}\binom{n-s}{q-r}.
\]
\end{enumerate}
\end{thm}
\begin{proof}
By Lemma \ref{lem:span}, the 
set of vectors $\cB^{q,r}_S(\bl)= \cV^{q,r}_S(\bl)
\bigcup \cW^{q,r}_S(\bl)$
spans the vector space $A^q(\cG^\ell_n)$.  So to establish this result, it 
suffices to show that these vectors are eigenvectors of the 
endomorphism $\omega^q_\bl(S,r)$, and that the dimensions of the eigenspaces 
are as asserted.  We will prove this by induction on $q-r$.  

For ease of notation, we will suppress dependence on $\bl$ 
in the proof, and, 
for instance, write simply $\omega^q(S,r)=\omega^q_\bl(S,r)$ 
and $\Psi^q_{S,r} = \Psi^q_{S,r}|_{y_j\mapsto\la_j}$.
Using Lemma~\ref{lem:sym action}, it suffices to consider the case $S \subset [n]$.  
Let $J\subset S$, $K \subset [n] \setminus S$, and recall that
\[ 
\begin{aligned}
\cV^{q,r}_{S} &= \{e_Je_K \mid |J| \le r-1\}
\bigcup \{\eta_{S} e_J e_K \mid|J| = r-1\} \quad \text{and} \\
\cW^{q,r}_{S} &= \{e_{S} e_K \ (\text{if}\ q \ge s)\} 
\bigcup \{(\partial{e_J}) e_K \mid |J| \ge r+1\} \bigcup
\{\eta_{S} e_J e_K \mid|J| \ge r\}. 
\end{aligned}
\]

In the case $q-r=0$, we have 
$\cV^{r,r}_{S} = \{e_Je_K \mid |J| \le r-1\}
\bigcup \{\eta_{S} e_J \mid|J| = r-1\}$, 
$\cW^{r,r}_S=\{\partial{e}_J \mid |J|=r+1\}$, and 
$\omega^r(S,r)=\Psi^r_{S,r}$.  By Lemma 
\ref{lem:Psi values}, if $|J|\le r-1$, then $\Psi^r_{S,r}(e_J e_K)=0$.  
If $|J|=r-1$, then, using Lemma \ref{lem:Psi values} again, we have
\[
\begin{aligned}
\Psi^r_{S,r}(\eta_S e_J) &= 
\sum_{i \in S} \la_i \Psi^r_{S,r}(e_i e_J)
=\sum_{i \in S} \la_i(\la_S e_i e_J - \eta_S \partial(e_i e_J)) \\
&= \la_S \eta_S e_J -\sum_{i\in S}\la_i \eta_S(e_J - e_i\partial e_J) \\
&= \la_S \eta_S e_J - \la_S \eta_S e_J + 
\eta_S \eta_S \partial e_J = 0.
\end{aligned}
\]
Thus, every element of
$E^{r}(0) = \Span \cV^{r,r}_S$ 
is a $0$-eigenvector of $\omega^r(S,r)$.  
A straightforward exercise reveals that $\dim E^{r}(0) = 
\sum_{k=0}^r\binom{s}{k}\binom{n-s}{r-k} - \binom{s-1}{r}$.  
If $|J|=r+1$, then, using Lemma \ref{lem:Psi values} again,
\[
\begin{aligned}
\Psi^r_{S,r}(\partial e_J) &= \sum_{k=1}^{r+1} (-1)^{k-1} 
\Psi^r_{S,r}(e_{J_k}) 
= \sum_{k=1}^{r+1} (-1)^{k-1}(\la_S e_{J_k} - \eta_S \partial e_{J_k}) \\
&=\la_S \partial e_J - \eta_S \partial^2 e_J = \la_S \partial e_J
\end{aligned}
\]
Thus, every element of 
$E^{r}(\la_S) = \Span \cW^{r,r}_S$ 
is a $\la_S$-eigenvector of $\omega^r(S,r)$.  Note that
$\dim E^{r}(\la_S)=\binom{s-1}{r}$.  Since $\dim E^{r}(0)+\dim E^{r}(\la_S)= \dim A^r(\cG^\ell_n)$, 
the above calculations establish Theorem \ref{thm:diag} in the case $q-r=0$.

If $q-r \ge 1$, then by induction, for each $k \ge 1$, $\omega^q(S,r+k)$
is diagonalizable, 
with eigenvalues $0$ and $\la_S$, and corresponding 
eigenspaces $E^{r+k}(0)=\Span \cV^{r+k,r}_S$ and 
$E^{r+k}(\la_S)=\Span \cW^{r+k,r}_S$.  
In the determination of the eigenstructure of $\omega^q(S,r)$, 
we will use the 
recursion provided by Lemma \ref{lem:recursion} in the following form:
\begin{equation} \label{eq:recursion}
\omega^q(S,r) = \Psi^q_{S,r} - \sum_{k=1}^{s-r-1} \binom{r+k-1}{k} 
\omega^q(S,r+k).
\end{equation}

First, consider the $0$-eigenspace of the endomorphism $\omega^q(S,r)$.  
If $|J| \le r-1$, then by \eqref{eq:recursion}, Lemma \ref{lem:Psi values}, 
and induction, we have
\[
\omega^q(S,r)(e_J e_K) = \Psi^q_{S,r}(e_J e_K) - \sum_{k=1}^{s-r-1} 
\binom{r+k-1}{k} \omega^q(S,r+k)(e_J e_K)=0.
\]
If $|J|=r-1$, then $\omega^q(S,r+k)(\eta_S e_J e_K)=0$ for $k \ge 1$ by 
Lemma \ref{lem:Psi values}.  Using \eqref{eq:recursion} 
and Lemma \ref{lem:Psi values}, we have
\[
\begin{aligned}
\omega^q(S,r)(\eta_S e_J e_K) &= \Psi^q_{S,r}(\eta_S e_J e_K) - 
\sum_{k=1}^{s-r-1} \binom{r+k-1}{k} \omega^q(S,r+k)(\eta_S e_J e_K)\\
&=\Psi^q_{S,r}(\eta_S e_J e_K) = 
\sum_{i\in S} \la_i \Psi^q_{S,r}(e_i e_J e_K)\\
&= \sum_{i\in S} \bigl[\la_i \la_S e_i e_J e_K - 
\la_i \eta_S \partial(e_i e_J) e_K\bigr]\\
&= \la_S \eta_S e_J e_K - \sum_{i\in S} \la_i \eta_S e_J e_K +
\sum_{i \in S} \la_i \eta_S e_i(\partial{e}_J)e_K\\
&= \la_S \eta_S e_J e_K - \la_S \eta_S e_J e_K + \eta_S \eta_S(\partial{e}_J)
e_K=0.
\end{aligned}
\]

Next, consider the $\la_S$-eigenspace.  If $q \ge s$, we must show that $e_S 
e_K$ is an eigenvector of $\omega^q(S,r)$ corresponding to the eigenvalue 
$\la_S$ 
for each $K \subset [n]\setminus S$ with $|K|=q-s$.  By induction, we have
$\omega^q(S,r+k)(e_S e_K) = \la_S e_S e_K$ for each $k \ge 1$.  By Lemma 
\ref{lem:Psi values}, we have 
$\Psi^q_{S,r}(e_S e_K) = 
\binom{s}{r} \la_S e_S e_K - 
\binom{s-1}{r-1} \eta_S \partial{e}_S e_K$.  Since 
$\eta_S \partial{e}_S = \la_S e_S$, we have $\Psi^q_{S,r}(e_S e_K) = 
\binom{s-1}{r}\la_S e_S e_K$.  Hence, by \eqref{eq:recursion}, we have
\[
\begin{aligned}
\omega^q(S,r)(e_S e_K)&= \Psi^q_{S,r}(e_S e_K) - \sum_{k=1}^{s-r-1} 
\binom{r+k-1}{k} \omega^q(S,r+k)(e_S e_K)\\
&=\binom{s-1}{r}\la_S e_S e_K - 
\sum_{k=1}^{s-r-1} \binom{r+k-1}{k} \la_S e_S e_K 
=\la_S e_S e_K,
\end{aligned}
\]
using the binomial identities
\begin{equation*} \label{eq:binom id}
\sum_{k=0}^{p} \binom{N+k}{k} = \binom{N+p+1}{p} = \binom{N+p+1}{N+1}.
\end{equation*}
with $N=r-1$ and $p=s-r-1$.

If $|J| \ge r+1$, we must show that $\omega^q(S,r)(\partial{e}_J e_K) =
\la_S \partial{e}_J e_K$.  Suppose $|J|=r+p+1$ for some $p\ge 0$.  Then, by 
Lemma \ref{lem:Psi values}, we have
\[
\begin{aligned}
\Psi^q_{S,r}(\partial{e}_J e_K)&=\sum_{i=1}^{r+p+1} (-1)^{i-1} 
\Psi^q_{S,r}(e_{J_i} e_K) \\
&=\sum_{i=1}^{r+p+1} (-1)^{i-1} \left[
\binom{r+p}{r} \la_S e_{J_i} e_K - \binom{r+p-1}{r-1} \eta_S 
(\partial{e}_{J_i}) e_K\right]\\
&= \binom{r+p}{r} \la_S (\partial{e}_J) e_K - \binom{r+p-1}{r-1} \eta_S 
(\partial^2{e}_{J}) e_K \\
&=\binom{r+p}{r} \la_S (\partial{e}_J) e_K.
\end{aligned}
\]
By induction, we have
\[
\omega^q(S,r+k)((\partial{e}_J)e_K)= 
\begin{cases}
\la_S (\partial{e}_J)e_K & \text{if $1 \le k \le p$},\\
0 & \text{if $p+1 \le k \le s-r-1$}.
\end{cases}
\]
So using the recursion \eqref{eq:recursion} and the identities \eqref{eq:binom 
id}, we obtain
\[
\begin{aligned}
\omega^q(S,r)((\partial{e}_J) e_K)&=
\binom{r+p}{r} \la_S (\partial{e}_J) e_K - \sum_{k=1}^p \binom{r+k-1}{k}
\la_S (\partial{e}_J) e_K \\
&= \la_S (\partial{e}_J) e_K.
\end{aligned}
\]

If $|J| \ge r$, we must show that $\omega^q(S,r)(\eta_S e_J e_K) =
\la_S \eta_S e_J e_K$.  Suppose $|J|=r+p$ for some $p\ge 0$.  Then, by 
Lemma \ref{lem:Psi values}, we have
\[
\begin{aligned}
\Psi^q_{S,r}(\eta_S e_J e_K)&=\sum_{i\in S} 
\la_i \Psi^q_{S,r}(e_i e_J e_K) \\
&=\sum_{i\in S} y_i \left[
\binom{r+p+1}{r} \la_S e_i e_J e_K - \binom{r+p}{r-1} \eta_S \partial(e_i e_J) 
e_K\right]\\
&= \binom{r+p+1}{r} \la_S \eta_S e_J e_K - \binom{r+p}{r-1}
\sum_{i\in S} \la_i \eta_S (e_J - e_i\partial{e}_J)e_K\\
&= \left[\binom{r+p+1}{r} 
-\binom{r+p}{r-1}\right] \la_S \eta_S e_J e_K+ \binom{r+p}{r-1} \eta_S \eta_S 
(\partial{e}_J) e_K\\
&= \binom{r+p}{r} \la_S \eta_S e_J e_K.
\end{aligned}
\]
By induction, we have
\[
\omega^q(S,r+k)(\eta_S e_J e_K)= 
\begin{cases}
\la_S \eta_S e_J e_K & \text{if $1 \le k \le p$},\\
0 & \text{if $p+1 \le k \le s-r-1$}.
\end{cases}
\]
So using the recursion \eqref{eq:recursion} and the identities \eqref{eq:binom 
id}, we obtain $\omega^q(S,r)(\eta_S e_J e_K)=\la_S \eta_S e_J e_K$ as above.

Thus the vectors 
in the sets $\cV^{q,r}_S(\bl)$ and $\cW^{q,r}_S(\bl)$ are 
eigenvectors of $\omega^q(S,r)$ 
corresponding to the eigenvalues $0$ and $\la_S$ as asserted.  Since these 
vectors span $A^q(\cG^\ell_n)$ by Lemma 
\ref{lem:span}, it remains to compute the dimensions of the eigenspaces 
$E^q(0)=\Span \cV^{q,r}_S$ and $E^q(\la_S)=\Span \cW^{q,r}_S$ corresponding to these eigenvalues.

If $|J|=p$ and $|J|+|K|=q$, then 
$\Span \{e_J e_K \mid J \subset S, K \subset [n]\setminus S\}$ 
has dimension 
\[
\binom{s}{p} \binom{n-s}{q-p}.
\]
If $|J|=p+1$ and $|J|-1+|K|=q$, then 
$\Span \{(\partial{e}_J) e_K \mid J \subset S, K \subset [n]\setminus S\}$ has dimension 
\[
\dim \im[\partial:A^{p+1}(\cG^s_s) \to A^p(\cG^s_s)]\cdot \binom{n-s}{q-p} =\binom{s-1}{p} \binom{n-s}{q-p}.
\]
If $|J|=p-1$ and $|J|+1+|K|=q$, then 
$\Span \{\eta_S e_J e_K \mid J \subset S, 
K \subset [n]\setminus S\}$ has dimension 
\[
\dim \ker[\eta_S:A^{p}(\cG^s_s) \to A^{p+1}(\cG^s_s)]\cdot \binom{n-s}{q-p} 
=\left[\binom{s}{p}-\binom{s-1}{p}\right] \binom{n-s}{q-p}.
\]

Using these calculations, it is readily checked that
\[
\begin{aligned}
\dim E^q(0) &= \sum_{p=0}^{r} \binom{s}{p}\binom{n-s}{q-p} -
\binom{s-1}{r}\binom{n-s}{q-r},\ \text{and} \\
\dim E^q(\la_S) &= \sum_{p=r+1}^{\min(q,s)} \binom{s}{p}\binom{n-s}{q-p} +
\binom{s-1}{r}\binom{n-s}{q-r}.
\end{aligned}
\]
The fact that $\dim E^q(0) + \dim E^q(\la_S) = 
\dim A^q(\cG^\ell_n) = \binom{n}{q}$ 
may be checked using the binomial identities
\[
\sum_{p=0}^k \binom{m}{p}\binom{N}{k-p} = \binom{m+N}{k} 
\ \text{and}\ 
\sum_{p=0}^m \binom{m}{p}\binom{N}{k+p} = \binom{m+N}{m+r} 
= \binom{m+N}{N-k}
\]
with $m=s$, $N=n-s$, and $k=q$ in the case $q<s$, and 
$m=s$, $N=n-s$, and $k=n-s-q$ in the case $q \ge s$.
\end{proof}

If $n=\ell$ and $\bl \neq 0$, the complex $(A^\bul(\cG),e_\bl)$ is acyclic.  
So assume that $n > \ell$.  Then, for $\bl \neq 0$, 
the cohomology of this complex is 
concentrated in dimension $\ell$, and $\dim H^\ell(\cG)=\binom{n-1}{\ell}$.  
Let $\rho=\rho_{_\cG}:A^\ell(\cG) \to H^\ell(\cG)$ denote the projection.  
Since $\omega^\bul_\bl(S,r)$ is a chain map,  the kernel of this projection, 
$\ker(\rho) \subset A^\ell(\cG)$, is an invariant subspace for $\omega^\ell_\bl(S,r)$.

\begin{lem} \label{lem:diag lem}
Let $T:V \to V$ be an endomorphism of a finite dimensional (complex) vector space, and $V'$ an invariant subspace.  If $T$ is diagonalizable, then the induced endomorphism $T''$ on the quotient $V''=V/V'$ is also diagonalizable, 
and the spectrum of $T''$ is contained in the spectrum of $T$.
\end{lem}
\begin{proof}
Let $T'$ denote the restriction of $T$ to $V'$, and let $\pi:V \to V''$ be the projection.  The vector space $V$ admits a basis $\cB=\{v_1,\dots,v_k,v_{k+1},\dots,v_n\}$ for which $\cB'=\{v_1,\dots,v_k\}$ 
is a basis for the subspace $V'$ and $\cB''=\{\pi(v_{k+1}),\dots,\pi(v_n)\}$ is a basis for the quotient $V''$.  The matrix of $T$ relative to the basis $\cB$ is 
\[
\sfA = \begin{pmatrix} \sfA' & *\phantom{''} \\ 
0\phantom{'} & \sfA'' \end{pmatrix},
\]
where $\sfA'$ is the matrix of $T'$ relative to $\cB'$ and $\sfA''$ is the matrix of the induced endomorphism $T''$ relative to $\cB''$.

Let $r_1,\dots,r_m$ be the distinct eigenvalues of $T$.  
Since $T$ is diagonalizable, the minimal polynomial $p$ of $T$ factors as
$p(t) = (t-r_1) \cdots (t-r_m)$.  The polynomial $p$ annihilates the matrix $\sfA$
of $T$, $p(\sfA)=0$.  Using the block decomposition of $\sfA$ above, it follows that $p$ also annihilates the matrix $\sfA''$ of $T''$, $p(\sfA'')=0$.  
Consequently, the minimal polynomial $p''$ of $T''$ divides $p$.  Hence, $p''$ is of the form $(t-r_{i_1}) \cdots (t-r_{i_j})$, $T''$ is diagonalizable, and
the eigenvalues of $T''$ are among the eigenvalues of $T$.
\end{proof}

For an arrangement $\A$ of arbitrary combinatorial type $\cT$, and 
$\cT$-nonresonant weights $\bl$, we recall the $\beta$\textbf{nbc} 
basis of \cite{FT} for the single nonvanishing cohomology group 
$H^\ell(\cT)=H^\ell(\sfM;\LL)$.  Recall that the hyperplanes of  
$\A=\{H_j\}_{j=1}^n$ are ordered.  
A circuit is an inclusion-minimal
dependent set of hyperplanes in $\A$, and a broken circuit is a set
$T$ for which there exists $H < \min(T)$ so that $T \cup\{H\}$ is a
circuit.  A frame is a maximal independent set, and an \textbf{nbc} frame
is a frame which contains no broken circuit.  Since $\A$ contains $\ell$
linearly independent hyperplanes, every frame has cardinality~$\ll$. 
The set of  \textbf{nbc} frames is a basis for $A^\ell(\cT)$.
An \textbf{nbc} frame
$B=(H_{j_1},\dots,H_{j_\ll})$ is a $\beta$\textbf{nbc} frame provided
that for each $k$, $1\le k \le \ll$, there exists $H\in \A$ such that
$H<H_{j_k}$ and $(B\setminus\{H_{j_k}\})\cup\{H\}$ is a frame.  
Note that these constructions depend only on the 
combinatorial type $\cT$ of $\A$, 
and let $\beta$\textbf{nbc}$(\cT)$ be the set of all $\beta$\textbf{nbc} 
frames of an arrangement of type $\cT$.  

\begin{definition}\label{def:bnbc}
Given $B=(H_{j_1},\dots,H_{j_\ll})$ in
$\beta$\rm{\textbf{nbc}}$(\cT)$, define $\xi(B) \in A^\ll(\cT)$ by
$\xi(B) = \wedge_{p=1}^\ll a_\la(X_p)$, where $X_p=\bigcap_{k=p}^\ll
H_{j_k}$ and $a_\la(X) = \sum_{X \subseteq H_i} \la_{i} a_{i}$.
Denote the cohomology class of $\xi(B)$ in
$H^\ell(\cT)=H^\ll(A^\bul(\cT),a_\bl)$ by the same symbol. 
The set $\{ \xi(B) \mid B\in 
\beta\text{\rm{\textbf{nbc}}}(\cT)\}$ is the $\beta$\rm{\textbf{nbc}} basis for $H^\ell(\cT)$.
\end{definition}

\begin{thm} \label{thm:Odiag} 
Let $S \subset [n+1]$ be a subset of cardinality $s$, and fix $r$,  $1\le r \le \min(\ell,s-1)$.  
For $\cG$-nonresonant weights $\bl$ satisfying $\la_S\neq 0$, the endomorphism
${\Omega}_\bl(S,r)$ of $H^\ell(\cG)$ induced by ${\omega}^\ell_\bl(S,r)$ 
is diagonalizable, with eigenvalues $0$ and $\la_S$.  The dimension of the $\la_S$-eigenspace is
\[
\sum_{p=r+1}^{\min(\ell,s)} \binom{s}{p} \binom{n-s-1}{\ell-p} + \binom{s-1}{r} \binom{n-s-1}{\ell-r},
\]
and the dimension of the $0$-eigenspace is
\[
\sum_{p=0}^r \binom{s}{p} \binom{n-s-1}{\ell-p} - \binom{s-1}{r} \binom{n-s-1}{\ell-r}.
\]
\end{thm}
\begin{proof}
By Theorem \ref{thm:diag} and Lemma \ref{lem:diag lem}, the endomorphism 
$\Omega_\bl(S,r)$ is diagonalizable, with spectrum contained in $\{0,\la_S\}$.

Let 
${\mathbf I}=\{I=(i_1,\ldots,i_\ell) \mid 1\leq i_1<i_2\cdots<i_\ell\leq n\}$.
Then $\{e_I\mid I\in \b{I}\}$ is the \textbf{nbc} basis of $A^\ell(\cG)$ 
and 
$\{\xi_I=\lambda_{i_1}\cdots \lambda_{i_\ell}e_I\mid I\in \b{I}, 1\notin I\}$ 
is the $\beta$\textbf{nbc} basis of $H^\ell(\cG)$.  
The projection $\rho:A^\ell(\cG) \to H^\ell(\cG)$ is given by
\[
\rho(e_I)=\begin{cases}
\phantom{-}(\la_{i_1} \cdots \la_{i_\ell})^{-1} \xi_I & \text{if $1 \notin I$,}\\
-(\la_{i_1} \cdots \la_{i_\ell})^{-1} \sum_{j \notin I} \xi_j \xi_{I_1} & \text{if $1 \in I$.}
\end{cases}
\]

Using Lemma \ref{lem:sym action}, we can assume that $S \subset [2,n]$.   
Since $\rho \circ \omega^\ell_\bl(S,r) = \Omega_\bl(S,r) \circ \rho$, if $v$ 
is an eigenvector of $\omega^\ell_\bl(S,r)$ and $\rho(v) \neq 0$, then $\rho(v)$ is an 
eigenvector of $\Omega_\bl(S,r)$.  
Let $J \subset S$ and $K \subset [2,n] \setminus S$.  Note that $1\notin K$.  Then one can check 
that the $0$-eigenspace of $\Omega_\bl(S,r)$ is spanned by 
\[
\{\rho(e_J e_K) \mid |J| \le r-1\} \bigcup 
\{\rho(\eta_S e_J e_K) \mid |J|=r-1\},
\]
that the $\la_S$-eigenspace of $\Omega_\bl(S,r)$ is spanned by 
\[
\{\rho(e_S e_K) \mid \text{if $\ell\ge s$}\} \bigcup 
\{\rho((\partial e_J) e_K) \mid |J| \ge r+1\} \bigcup
\{\rho(\eta_S e_J e_K) \mid |J| \ge r\},
\]
and that the dimensions of these eigenspaces are as asserted.
\end{proof}

\begin{exm} \label{exm:345}
Let $n=5$, $\ell=2$, $S=\{3,4,5\}$, and $r=1$.  
By Theorem \ref{thm:Odiag}, 
for $\cG$-nonresonant weights satisfying $\la_S \neq 0$, 
the endomorphism 
$\Omega_\bl(S,r)$ of $H^2(\cG) \simeq \C^6$ is diagonalizable, 
the $\la_S$-eigenspace 
is $5$-dimensional, and the $0$-eigenspace is $1$-dimensional (note that $\binom{p}{q}=0$ if $p<q$).  Calculating as in the proof of Theorem~\ref{thm:Odiag}, we find that the $\la_S$-eigenspace has basis
\begin{alignat*}{3}
\rho(\la_2\la_3\la_5(\partial{e}_{3,5})e_2)&=\la_5 \xi_{2,3}-\la_3\xi_{2,5},
&\rho(-\la_3 \eta_{3,4,5} e_3) &= \xi_{3,4}+\xi_{3,5},
\\
\rho(\la_2\la_4\la_5(\partial{e}_{4,5})e_2)&=\la_5 \xi_{2,4}-\la_4\xi_{2,5}, 
&\rho(\la_5 \eta_{3,4,5} e_5) &= \xi_{3,5}+\xi_{4,5}, \\
\rho(\la_3\la_4\la_5\partial{e}_{3,4,5})&=
\la_5\xi_{3,4}-\la_4\xi_{3,5}+\la_3\xi_{4,5},\quad &
\end{alignat*}
and the $0$-eigenspace has basis 
\[
\rho(\la_1\la_2 e_{1,2})=\xi_{2,3}+\xi_{2,4}+\xi_{2,5}.
\]
\end{exm}

\section{Nonresonant eigenvalues} \label{sec:evalues}

In this section, we prove that the Gauss-Manin endomorphism 
$\GM_\bl(\cT',\cT)$ of \eqref{eq:1form} is diagonalizable and determine 
its eigenvalues.
We accomplish this by showing  that 
the endomorphism  $\widetilde{\Omega}_\bl(\cT',\cT)$ in
the commutative diagram 
\eqref{eq:main1} may be replaced by the endomorphism 
${\Omega}_\bl(S,r)$, whose eigenstructure was computed in Theorem~\ref{thm:Odiag}.

For an arbitrary type $\cT$, 
let $I^\bullet(\cT)$ be the corresponding Orlik-Solomon 
ideal, so that $A^\bullet(\cT) \simeq A^\bullet(\cG) / I^\bullet(\cT)$.  
The natural projection of $A^\bullet(\cG)$ onto $A^\bullet(\cT)$ is a chain 
map $\pi:(A^\bullet(\cG),e_\bl) \to (A^\bullet(\cT),a_\bl)$ which,
for $\cT$-nonresonant weights $\bl$, induces 
the projection $\tau:H^\ell(\cG) \to H^\ell(\cT)$ upon passage to cohomology.  
If ${\rho_{_\cG}}:A^\ell(\cG) \to H^\ell(\cG)$ and ${\rho_{_\cT}}:A^\ell(\cT) \to H^\ell(\cT)$ 
are the projections, then $\tau \circ \rho_{_\cG} = \rho_{_\cT} \circ \pi$.

\begin{thm}
\label{thm:decomp}
If $\cT'$ is a degeneration of $\cT$ with principal dependence $(S,r)$, then 
$\Omega_\bl(\cT',\cT) \circ \tau = \tau \circ \Omega_\bl(S,r)$.  
In other words, the following diagram commutes:
\[
\begin{CD}
H^\ell(\cG) @>\tau>> H^\ell(\cT)\\
@VV{\Omega}_\bl(S,r)V      
@VV\Omega_\bl(\cT',\cT)V \\
H^\ell(\cG) @>\tau>> H^\ell(\cT)
\end{CD}
\]
\end{thm}
\begin{proof}
As noted in the introduction, the Gauss-Manin endomorphism 
$\Omega_\bl(\cT',\cT)$ of $H^\ell(\cT)$ is induced by the endomorphism 
$\widetilde\Omega_\bl(\cT',\cT)$ of $H^\ell(\cG)$, see \cite[Thm.~7.3]{CO4} and \eqref{eq:main1}.  
In turn, $\widetilde\Omega_\bl(\cT',\cT)$ is the map in cohomology induced by 
the cochain endomorphism $\omega_\bl^\bul(\cT',\cT)$ 
of the complex $(A^\bullet(\cG),e_\bl)$, see Theorem \ref{thm:induced endo}.
The map $\omega_\bl^\bul(\cT',\cT)$ also induces a cochain endomorphism 
$\bar\omega_\bl^\bul(\cT',\cT)$ of $(A^\bullet(\cT),a_\bl)$, and the Gauss-Manin endomorphism 
$\Omega_\bl(\cT',\cT)$ may be realized as the map in cohomology induced by the latter, 
see \cite[Thm.~7.1]{GM3}.  In summary, we have the following 
commutative diagram.  
\begin{equation} \label{eq:fancy}
\begin{diagram}[notextflow]
   A^\ell(\cG)  &    &\rTo^{\pi} &      &   A^\ell(\cT)   \\
      & \rdTo_{\rho_{_\cG}} &      &      & \vLine^{\bar\omega_\bl^\ell(\cT',\cT)}& \rdTo_{\rho_{_\cT}}  \\
\dTo^{\omega^\ell_\bl(\cT',\cT)} &    &   H^\ell(\cG)   & \rTo^\tau  & \HonV   &    &  H^\ell(\cT)  \\
      &    & \dTo^{\widetilde\Omega_\bl(\cT',\cT)}  &      & \dTo   \\
    A^\ell(\cG)  & \hLine & \VonH   & \rTo^{\pi} &    A^\ell(\cT)   &    & \dTo_{\Omega_\bl(\cT',\cT)} \\
      & \rdTo_{\rho_{_\cG}} &      &      &      & \rdTo_{\rho_{_\cT}}  \\
      &    &   H^\ell(\cG)   &      & \rTo^\tau  &    &  H^\ell(\cT)  \\
\end{diagram}
\end{equation}
To establish the theorem, it suffices to show that the endomorphisms $\omega_\bl^\bul(\cT',\cT)$ and 
$\omega_\bl^\bul(S,r)$ of $A^\bul(\cG)$ induce the same endomorphism of $A^\bul(\cT)$.

The Orlik-Solomon ideal $I^\bul(\cT)$ gives rise to a subcomplex 
$I^\bullet_R(\cT) = I^\bullet(\cT) \otimes R$ of the 
Aomoto complex $A^\bullet_R(\cG)$, with quotient $A^\bullet_R(\cT)$, 
the Aomoto complex of type $\cT$.  Since  
$\omega_\bl^\bul(\cT',\cT)$ and $\omega_\bl^\bul(S,r)$ are specializations at $\bl$ of 
the corresponding endomorphims of the Aomoto complex $A^\bul_R(\cG)$, it is enough 
to show that $\omega^\bul(\cT',\cT)$ and 
$\omega^\bul(S,r)$ induce the same endomorphism of $A_R^\bul(\cT)$.

By Theorem \ref{thm:principal}, there are dependence pairs $(S_i,r_i)$, 
$1\le i \le k$, such that $\Dep(\cT(S_i,r_i))^* \subset \Dep(\cT)^*$ and 
$\Dep(\cT')^*= 
\bigcup_{i=0}^k \Dep(\cT(S_i,r_i))^*$, where $(S_0,r_0)=(S,r)$ is the pair of 
principal dependence.  
It follows that there are constants $c_i$ so that
${\omega}^\bullet(\cT') = {\omega}^\bullet(S,r) + \sum_{i=1}^k c_i \cdot 
{\omega}^\bullet(S_i,r_i)$.

If $\Dep(\cT(S_i,r_i))^* \subset \Dep(\cT)^*$, it follows from  
Theorem \ref{thm:diag} that the image of 
${\omega}^\bullet(S_i,r_i):A^\bullet_R(\cG)\to A^\bullet_R(\cG)$ 
is contained in $I^\bullet_R(\cT)$.  Consequently, the endomorphisms 
$\bar\omega^\bullet(\cT')$ 
and $\bar\omega^\bullet(S,r)$ of the Aomoto complex $A^\bullet_R(\cT)$ induced 
by 
${\omega}^\bullet(\cT')$ and ${\omega}^\bullet(S,r)$ are equal.

Finally,
${\omega}^\bullet(\cT')=
{\omega}^\bullet(\cT',\cT)+{\omega}^\bullet(\cT)$, 
see \eqref{eq:GP endo}. It follows from the definitions 
that the image of ${\omega}^\bullet(\cT)$ is also
contained in 
$I^\bullet_R(\cT)$.  
Hence, the endomorphisms $\bar\omega^\bullet(\cT')$ 
and $\bar\omega^\bullet(\cT',\cT)$ of $A^\bullet_R(\cT)$ induced by 
${\omega}^\bullet(\cT')$ and ${\omega}^\bullet(\cT',\cT)$ are equal.
\end{proof}

Theorem \ref{thm:Odiag} and Theorem \ref{thm:decomp} 
yield the result stated in the introduction.

\begin{thm}\label{thm:evalues} 
Let $\cT'$ be a degeneration of $\cT$ with principal dependence $(S,r)$, and 
$\bl$ a collection of $\cT$-nonresonant weights satisfying $\la_S \neq 0$.  
Then the Gauss-Manin endomorphism $\Omega_\bl(\cT',\cT)$ is diagonalizable, with spectrum contained in $\{0,\la_S\}$.
\end{thm} 
\begin{proof}
By Theorem \ref{thm:Odiag}, the endomorphism $\Omega_\bl(S,r)$ of $H^\ell(\cG)$ 
is diagonalizable, with eigenvalues $0$ and $\la_S$.  By Theorem \ref{thm:decomp}, 
we have $\Omega_\bl(\cT',\cT) \circ \tau = \tau \circ \Omega_\bl(S,r)$.  
Checking that $\ker(\tau) \subset H^\ell(\cG)$ is an invariant subspace 
for $\Omega_\bl(S,r)$, the result follows from Lemma \ref{lem:diag lem}.
\end{proof}

\begin{rem}\label{rem:algorithm}
The Gauss-Manin endomorphism $\Omega_\bl(\cT',\cT)$ of $H^\ell(\cT)$ 
is determined by the endomorphism $\Omega_\bl(S,r)$ of $H^\ell(\cG)$ 
and the projection $\tau:H^\ell(\cG) \to H^\ell(\cT)$ via the equality  
$\Omega_\bl(\cT',\cT) \circ \tau = \tau \circ \Omega_\bl(S,r)$.  
Together 
with the explicit description of the eigenstructure of $\Omega_\bl(S,r)$ provided by 
Theorems \ref{thm:diag} and \ref{thm:Odiag}, this yields an algorithm for 
finding the (geometric) multiplicities of the eigenvalues of $\Omega_\bl(\cT',\cT)$.
\end{rem}

The Gauss-Manin connection $\nabla = \sum \Theta_{\cT'} \otimes \Omega_\bl(\cT',\cT)$ 
on the vector bundle $\b{H} \to \sfB(\cT)$ with fiber $H^\ell(\cT)$ corresponds to a monodromy 
representation $\Psi:\pi_1(\sfB(\cT)) \to \Aut_\C\bigl(H^\ell(\cT)\bigr)$.  For a degeneration 
$\cT'$ of $\cT$, let $\gamma_{\cT'} \in \pi_1(\sfB(\cT))$ be a simple loop in $\sfB(\cT)$ 
around a generic point in $\sfB(\cT')$.  Then the automorphism $\Psi(\gamma_{\cT'})$ 
is conjugate to $\exp\bigl(-2 \pi \ii \Omega_\bl(\cT',\cT)\bigr)$, see for instance \cite[Prop.~4.1]{CO3}.  
Theorem~\ref{thm:evalues} yields:

\begin{cor}
Let $\cT'$ be a degeneration of $\cT$ with principal dependence $(S,r)$, and 
$\bl$ a collection of $\cT$-nonresonant weights satisfying $\la_S \neq 0$.  
Then the automorphism $\Psi(\gamma_{\cT'})$ is diagonalizable, 
with spectrum contained in $\{1,\exp(-2\pi \ii \la_S)\}$.
\end{cor}

We conclude with several examples which illustrate these results.

\subsection{Codimension zero} \label{subsec:codim zero}
Recall that $\cG$ denotes the combinatorial type of a general position arrangement 
of $n$ hyperplanes in $\C^\ell$, and that $n\ge \ell$.  
Weights $\bl=(\la_1,\dots,\la_n)$ are $\cG$-nonresonant if $\la_j\neq 0$ for each $j$.  
If $n = \ell$, then $H^\bul(\cG)=0$,  so we assume that $n>\ell$.  
Then $\dim H^\ell(\cG)=\binom{n-1}{\ell}$.  
The moduli space $\sfB(\cG)$ 
has codimension zero in $(\CP^\ell)^n$, and 
consists of all matrices $\sfb$ for which every $(\ell+1)\times(\ell+1)$ minor is nonzero, see
\eqref{eq:point}.  For general position arrangements, the Gauss-Manin connection was determined 
by Aomoto and Kita \cite{AK}.  The corresponding connection $1$-form is given by
$\nabla = \sum \Theta_\cT \otimes \Omega_\bl(\cT,\cG)$,
where the sum is over all $\ell+1$ element subsets $S$ of $[n+1]$, 
$\cT=\cT(S,\ell+1)$, and 
$\Theta_\cT$ is a logarithmic $1$-form on $(\CP^\ell)^n$ with a simple pole 
along the divisor defined by the vanishing of the $(\ell+1)\times(\ell+1)$ minor of $\sfb$ with rows indexed by $S$.  Theorem \ref{thm:Odiag} gives:

\begin{prop} \label{prop:GPdiag} 
Let $S$ be an $\ell+1$ element subset of $[n]$, 
let $\cT=\cT(S,\ell+1)$, 
and $\bl$ a collection of  
$\cG$-nonresonant weights satisfying $\la_S\neq 0$.  
Then the Gauss-Manin endomorphism $\Omega_\bl(\cT,\cG)$ 
is diagonalizable, with eigenvalues $0$ and $\la_S$.  
The dimension of the $\la_S$-eigenspace is $1$, 
and the dimension of the $0$-eigenspace is $\binom{n-1}{\ell}-1$.
\end{prop}

\subsection{Codimension one} \label{subsec:codim one}
If $\cT$ is a combinatorial type for which the cardinality of $\Dep(\cT)_{\ell+1}$ is $1$, 
then the moduli space $\sfB(\cT)$ is of codimension one in $(\CP^\ell)^n$.  
Write $\Dep(\cT)_{\ell+1} = \{K\}$.  As shown by Terao \cite{T1}, noted in the proof of 
Theorem~\ref{thm:principal}, and illustrated in Example \ref{ex:codim1}, 
the combinatorial type $\cT$ admits three types 
of degeneration $\cT'=\cT(S,r)$.  The principal dependencies of these degenerations are as follows.
\begin{enumerate}
\item[I:] $(S,\ell)$, where $|S|=\ell+1$ and $|S\cap K| \le \ell-1$;
\item[II:] $(S,\ell-1)$, where $S=K_p$, for each $p$, $1\le p \le \ell+1$;
\item[III:] $(S,\ell)$, where $S=(m,K)$, for each $m \in [n+1]\setminus K$.
\end{enumerate}

For the combinatorial type $\cT$ and $\cT$-nonresonant weights $\bl$, the 
Gauss-Manin connection was determined by Terao \cite{T1}.  
The corresponding connection $1$-form is given by 
$\nabla = \sum \Theta_{\cT'} \otimes \Omega_\bl(\cT',\cT)$, 
where $\cT'$ ranges over the three types of degeneration of $\cT$ noted above.  
In \cite{T1}, Terao also found the eigenvalues of the endomorphism $\Omega_\bl(\cT',\cT)$ 
and their algebraic multiplicities.  If $\bl$ 
satisfies $\la_S\neq 0$ 
for each of the principal dependence sets $S$ recorded above, Terao's result concerning the 
eigenstructure of the endomorphism $\Omega_\bl(\cT',\cT)$ may be 
strengthened as follows.

\begin{prop}
Let $\cT$ be a combinatorial type of codimension one,  
let $\cT'=\cT(S,r)$ be a degeneration of $\cT$, 
and $\bl$ a collection of  $\cT$-nonresonant weights satisfying $\la_S\neq 0$.  
Then the Gauss-Manin endomorphism $\Omega_\bl(\cT',\cT)$ is diagonalizable, with eigenvalues 
$0$ and $\la_S$.
\begin{enumerate}
\item If $\cT'$ is a degeneration of type {\rm{I}}, the dimension of the $\la_S$-eigenspace is $1$, and the 
dimension of the $0$-eigenspace is $\dim H^\ell(\cT)-1 = \binom{n-1}{\ell}-2$.
\item If $\cT'$ is a degeneration of type {\rm{II}}, the dimension of the $\la_S$-eigenspace is $n-\ell-1$, 
and the dimension of the $0$-eigenspace is $\binom{n-1}{\ell}-n+\ell$.
\item If $\cT'$ is a degeneration of type {\rm{III}}, the dimension of the $\la_S$-eigenspace is $\ell$, 
and the dimension of the $0$-eigenspace is $\binom{n-1}{\ell}-\ell-1$.
\end{enumerate}
\end{prop}
\begin{proof}
By Theorem \ref{thm:evalues}, the endomorphism $\Omega_\bl(\cT',\cT)$ is diagonalizable, with spectrum contained in $\{0,\la_S\}$.

Without loss, assume that $\Dep(\cT)_{\ell+1} = \{K\}$, where $K=[\ell+1]$.  Then the \textbf{nbc} basis of $A^\ell(\cT)$ consists of monomials $a_I$, where $I \subset [n]$, $|I|=\ell$, and $I \neq [2,\ell+1]$.  
Write $F=[2,\ell+1]$.  
The projection $\pi:A^\ell(\cG) \to A^\ell(\cT)$ is given by 
\[
\pi(e_I)=\begin{cases}
a_I & \text{if $I \neq F$,}\\
a_1 \partial a_F & \text{if $I=F$.}
\end{cases}
\]
The $\beta$\textbf{nbc} basis for $H^\ell(\cT)$ consists of monomials $\xi_I$, where 
$I \subset [2,n]$, $|I|=\ell$, and $I \neq F$.  The projection $\rho=\rho_{_\cT}:A^\ell(\cT) \to
H^\ell(\cT)$ is given by
\[
\rho(a_I)=\begin{cases}
\phantom{-}(\la_{i_1}\cdots \la_{i_\ell})^{-1} \xi_I & \text{if $1\notin I$,}\\
-(\la_{i_1}\cdots \la_{i_\ell})^{-1}\sum_{j \notin I} \xi_j \xi_{I_1} &
\text{if $1 \in I, I\not\subset K$,}\\
-(\la_K \la_{i_1}\cdots \la_{i_\ell})^{-1}
\sum_{j\notin K} \Bigl[\la_I\xi_j \xi_{I_1}
+ \xi_j \xi_p \partial\xi_{I_1}
\Bigr]&
\text{if $1 \in I, I= K\setminus\{p\}$.}
\end{cases}
\]

If $\cT'$ is a degeneration of type I with principal dependence $(S,\ell)$, 
then $|S \cap K| \le \ell-1$ and we can assume that $S \cap K \subset [3,\ell+1]$.  
By Theorem \ref{thm:diag}, the endomorphism $\omega_\bl^\ell(S,\ell)$ 
of $A^\ell(\cG)$ is diagonalizable, 
with eigenvalues $0$ and $\la_S$.  The  
$0$-eigenspace is spanned by $\{e_J e_L \mid |J| \le \ell-1\} \bigcup \{\eta_S e_J
\mid |J|=\ell-1\}$, where $J\subset S$ and $L\subset [n] \setminus S$, and the 
$\la_S$-eigenspace is spanned by $\partial{e}_S$.  By Theorem \ref{thm:decomp}, 
the endomorphism $\Omega_\bl(\cT',\cT)$ of $H^\ell(\cT)$ satisfies 
$\Omega_\bl(\cT',\cT) \circ \rho \circ \pi = \rho \circ \pi \circ \omega_\bl^\ell(S,\ell)$, 
see \eqref{eq:fancy}.
Write $S=(s_1,\dots,s_{\ell+1})$.   Calculations with the projections $\pi$ and $\rho$
yield
\[
\begin{aligned}
\rho \circ \pi (\partial{e}_S) &= (\la_{s_1} \cdots \la_{s_{\ell+1}})^{-1} \partial\xi_S,\\
\rho \circ \pi (e_J e_L) &= (\la_{i_1} \cdots \la_{i_{\ell}})^{-1} \xi_J \xi_L,
\ \text{where $I=(J,L)$ and $1\notin L$,}\\
\rho \circ \pi (\eta_S e_J) &= \la_{s_1}\la_{s_p}
(\la_{s_1} \cdots \la_{s_{\ell+1}})^{-1}(\xi_{S_p} \pm \xi_{S_1}),
\ \text{where $J=(s_2,\dots,\hat{s}_p,\dots,s_{\ell+1})$.}
\end{aligned}
\]
Checking that
\[
\{\partial\xi_S)\}\bigcup
\{\xi_J \xi_L \mid J \subset S, |J|\le \ell-1,L \subset [2,n]\setminus S\}
\bigcup \{\xi_{S_p} \pm \xi_{S_1} \mid 2\le p \le \ell+1 \}
\]
forms a basis for $H^\ell(\cT)$, we conclude 
that the dimensions of the eigenspaces are as asserted for a degeneration 
of type I.

If $\cT'$ is a degeneration of type II with principal dependence $(S,\ell-1)$, 
we can assume that $S=K_1=[2,\ell+1]$.  
By Theorem \ref{thm:diag}, the endomorphism $\omega_\bl^\ell(S,\ell-1)$ 
of $A^\ell(\cG)$ is diagonalizable, 
with eigenvalues $0$ and $\la_S$.  The  
$0$-eigenspace is spanned by $\{e_J e_L \mid |J| \le \ell-2\} \bigcup \{\eta_S e_J e_q
\mid |J|=\ell-2\}$, where $J\subset S$, $L\subset [n] \setminus S$, $q \notin S$, and the 
$\la_S$-eigenspace is spanned by $\{e_S\} \bigcup \{(\partial{e}_S)e_q \mid q \notin S\}$.
Note that the $\la_S$-eigenspace of $\omega_\bl^\ell(S,\ell-1)$ has dimension $n-\ell+1$.  
Note also that the $\la_S$-eigenvectors $\partial{e}_K=e_S-e_1 \partial{e}_S$ 
and $e_\bl \partial{e}_S$ are annihilated by the projection $\rho \circ \pi$.  
On the other hand, it is readily checked that  
\begin{equation} \label{eq:II S}
\{\rho\circ \pi ((\partial{e}_S) e_q) \mid \ell+2\le q \le n\}
\end{equation} 
is a linearly independent set of $(n-\ell-1)$ $\la_S$-eigenvectors for 
$\Omega_\bl(\cT',\cT)$ in $H^\ell(\cT)$.  Additionally, one can check that the set
\begin{equation} \label{eq:II 0}
\{\rho\circ\pi(e_J e_L) \mid J \subset S, |J|\le \ell-2\} \bigcup
\{\rho\circ\pi(\eta_S e_J e_q) \mid J \subset [3,\ell+1], |J|=\ell-2\}
\end{equation}
where $L \subset [n] \setminus S$ and $q\notin S$, is a linearly independent set of 
$0$-eigenvectors for $\Omega_\bl(\cT',\cT)$ in $H^\ell(\cT)$.  Checking that the dimension 
of the subspace spanned by the vectors \eqref{eq:II 0} is $\dim H^\ell(\cT)-(n-\ell-1)$, since 
eigenvectors associated to distinct eigenvalues are linearly independent, 
the vectors \eqref{eq:II S} and \eqref{eq:II 0} form a basis for $H^\ell(\cT)$.  Hence, 
the dimensions of the eigenspaces are as asserted for a degeneration of type II.

If $\cT'$ is a degeneration of type III with principal dependence $(S,\ell)$, 
we can assume that $S=K\bigcup\{q\}$ for some $q \in [\ell+2,n]$.  
By Theorem \ref{thm:diag}, the endomorphism $\omega_\bl^\ell(S,\ell)$ 
of $A^\ell(\cG)$ is diagonalizable, 
with eigenvalues $0$ and $\la_S$.  The  
$0$-eigenspace is spanned by $\{e_J e_L \mid |J| \le \ell-1\} \bigcup \{\eta_S e_J
\mid |J|=\ell-2\}$, where $J\subset S$, $L\subset [n] \setminus S$, and the 
$\la_S$-eigenspace is spanned by $\{\partial{e}_J \mid J \subset S, |J|=\ell+1\}$.
Note that the $\la_S$-eigenspace of $\omega_\bl^\ell(S,\ell)$ has dimension $\ell+1$.  
Note also that the $\la_S$-eigenvector $\partial{e}_K$ 
is annihilated by the projection $\rho \circ \pi$.  
Recall that $F=[2,\ell+1]$.  
Let $S_q$ denote the subspace of $H^\ell(\cT)$ 
spanned by $\{\xi_I \mid I \subset F\bigcup\{q\}\}$, and let ${\sf{p}}_q:H^\ell(\cT)
\to S_q$ be the natural projection.  For $J\subset F$, $|J|=\ell-1$, a calculation reveals that
${\sf{p}}_q \circ \rho \circ \pi(\eta_K e_J e_q) = \la_S (\la_2\cdots \la_{\ell+1} \la_q)^{-1}\xi_J \xi_q$.
Consequently, the set $\{\rho \circ \pi(\eta_K e_J e_q) \mid J\subset F, |J|=\ell-1\}$ 
is a linearly independent set of $\ell$ $\la_S$-eigenvectors for 
$\Omega_\bl(\cT',\cT)$ in $H^\ell(\cT)$.
Check that the set 
$\{\rho \circ \pi (e_J e_L) \mid J \subset S_1, |J|\le\ell-1,L\subset [n]\setminus S\}$ 
is a linearly independent set of $\dim H^\ell(\cT)-\ell$ $0$-eigenvectors for 
$\Omega_\bl(\cT',\cT)$ in $H^\ell(\cT)$.  It follows that 
the dimensions of the eigenspaces are as asserted for a degeneration of type III.
\end{proof}

\subsection{Further examples} \label{subsec:ex}
We present three examples of higher codimension. 

\begin{exm} \label{exm:Selberg}
Let $\cS$ be the
combinatorial type of the Selberg arrangement $\A$ in $\C^2$ with
defining polynomial
$
Q(\A)=u_1u_2(u_1-1)(u_2-1)(u_1-u_2)
$
depicted in
Figure \ref{fig:selberg}.
See \cite{Ao,SV,JK} for detailed studies of the Gauss-Manin
connections arising in the context of Selberg arrangements.

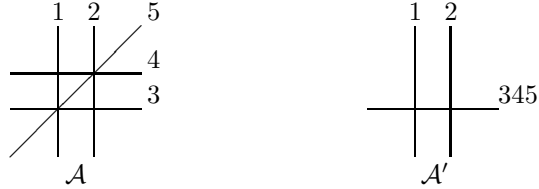
\begin{figure}[h]
\setlength{\unitlength}{.45pt}
\begin{picture}(300,130)(-200,-110)
\put(-200,-100){\line(0,1){110}}\put(-270,-100){\line(1,1){110}}
\put(-230,-100){\line(0,1){110}}\put(-270,-30){\line(1,0){110}}
\put(-270,-60){\line(1,0){110}}
\put(-155,15){5}\put(-155,-25){4}\put(-155,-55){3}\put(-235,15){1}
\put(-205,15){2}
\put(-224,-122){${\mathcal A}$}

\put(100,-100){\line(0,1){110}}
\put(70,-100){\line(0,1){110}}
\put(30,-60){\line(1,0){110}}
\put(140,-55){345}\put(65,15){1}\put(95,15){2}
\put(75,-122){${\mathcal A}'$}
\end{picture}
\caption{A Selberg arrangement and one degeneration}
\label{fig:selberg}
\end{figure} 

Here $\Dep(\cS)^*=\{126,346,135,245\}$.  Weights $\bl$ are $\cS$-nonresonant if
\[
\la_j\ (1\le j\le 6),\ \la_1+\la_2+\la_6,\ \la_1+\la_3+\la_5,
\ \la_2+\la_4+\la_5,\ \la_3+\la_4+\la_6\notin\Z_{\ge 0}.
\]
For 
$\cS$-nonresonant weights,  the 
$\beta\textbf{nbc}$ basis for $H^2(\cS)$ is 
$\{\Xi_{2,4},\Xi_{2,5}\}$, 
where $\Xi_{2,j}=(\la_2 a_2 +\la_4 a_4+\la_5 a_5)\la_j a_j$, 
see Definition \ref{def:bnbc}.
Recall that $\la_J=\sum_{j\in J}\la_j$.
The  projection map  
$\tau:H^2(\cG) \rightarrow H^2(\cS)$ is given by
\[
\tau(\xi_{i,j})=
\begin{cases}
-\Xi_{2,4}-\Xi_{2,5} &\text{if 
$(i,j)=(2,3)$,}\\
(\la_{2,4}\Xi_{2,4}+\la_4\Xi_{2,5})/\la_{2,4,5} & \text{if $(i,j)=(2,4)$,}\\
(\la_5\Xi_{2,4}+\la_{2,5}\Xi_{2,5})/\la_{2,4,5} &\text{if 
$(i,j)=(2,5)$,}\\
0 &\text{if $(i,j)=(3,4)$,}\\
(-\la_{5}\Xi_{2,4}-\la_{3,5}\Xi_{2,5})/\la_{1,3,5}&\text{if 
$(i,j)=(3,5)$,}\\
(-\la_{5}\Xi_{2,4}+\la_{4}\Xi_{2,5})/\la_{2,4,5}&\text{if 
$(i,j)=(4,5)$.}
\end{cases}
\]

The arrangement $\A'$ in Figure \ref{fig:selberg} represents one 
degeneration type
$\cS'$ of $\cS$.  Here
$\Dep(\cS',\cS)^*=\{34,35,45,134,145,234,235,345,356,456\}$.  
The sets $34$, $35$, $45$, and $345$ have $r=1$, and the others $r=2$.
The principal dependence is $(S,r)$, where $S=345$ and $r=1$.  
For $\cS$-nonresonant weights with $\la_S \neq 0$,
$\GM_\bl(\cS',\cS)$ is diagonalizable, 
with spectrum contained in $\{0,\la_S\}$ by Theorem \ref{thm:evalues}.  
The projection $\tau$ annihilates the $0$-eigenspace of $\Omega_\bl(S,r)$, 
and restricts to a surjection $E(\la_S) \twoheadrightarrow H^2(\cS)$, where 
$E(\la_S)$ is the $\la_S$-eigenspace of $\Omega_\bl(S,r)$, see Example \ref{exm:345}.  It follows that $\Omega_\bl(\cS',\cS)$ has eigenvalues $\la_S,\la_S$.  Note that $0$ is not an eigenvalue of
$\GM_\bl(\cS',\cS)$ in this instance.
\end{exm}

Although the eigenvalues are determined by the
principal dependence $(S,r)$,
the same  principal dependence may occur
for degenerations of different types. Thus the multiplicities of the
eigenvalues depend on the combinatorial
types as well. 

\begin{exm} \label{exm:not Selberg}
Consider the arrangement $\bar\A$ of type 
$\cT$ obtained from the arrangement $\A$ in Example \ref{exm:Selberg} 
by rotating line 1 by a (small)
angle about the triple point $135$, see Figure \ref{fig:not selberg}. 
Here, 
lines $1$ 
and $2$ meet in affine space, so $126$ is no longer dependent.
This change implies that $\dim A^2(\cT)=7$ and $\dim H^2(\cT)=3$. 

\begin{figure}[h]
\setlength{\unitlength}{.45pt}
\begin{picture}(300,130)(-200,-110)
\put(-200,10){\line(0,-1){115}}
\put(-275,-105){\line(1,1){112}}
\put(-190,-100){\line(-1,1){85}}
\put(-275,-30){\line(1,0){115}}
\put(-275,-60){\line(1,0){115}}
\put(-155,15){5}\put(-155,-25){4}\put(-155,-55){3}
\put(-280,-10){1}
\put(-205,15){2}
\put(-224,-122){${\bar\A}$}
\put(110,-100){\line(-1,1){85}}
\put(100,-105){\line(0,1){110}}
\put(30,-60){\line(1,0){110}}
\put(140,-55){345}\put(20,-10){1}\put(95,15){2}
\put(75,-122){${\bar\A}'$}
\end{picture}
\caption{A line arrangement and one degeneration}
\label{fig:not selberg}
\end{figure}
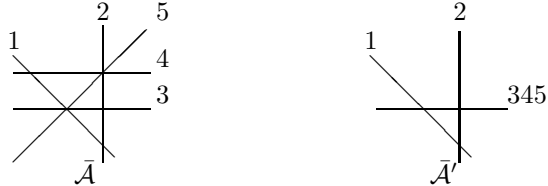 

Weights $\bl$ are $\cT$-nonresonant if
\[
\la_j\ (1\le j\le 6),\ \la_1+\la_3+\la_5,
\ \la_2+\la_4+\la_5,\ \la_3+\la_4+\la_6\notin\Z_{\ge 0}.
\]
For 
$\cT$-nonresonant weights, 
the $\beta\textbf{nbc}$ basis for $H^2(\cT)$ is $\{\Xi_{2,3},\Xi_{2,4},\Xi_{2,5}\}$, where 
$\Xi_{2,3}=\la_2 \la_3 a_{2,3}$ and 
$\Xi_{2,j}=(\la_2 a_2 +\la_4 a_4+\la_5 a_5)\la_j a_j$ for $j=4,5$.
The  projection map  
$\tau:H^2(\cG) \rightarrow H^2(\cT)$ is given by
\[
\tau(\xi_{i,j})=
\begin{cases}
\Xi_{2,3} &\text{if 
$(i,j)=(2,3)$,}\\
(\la_{2,4}\Xi_{2,4}+\la_4\Xi_{2,5})/\la_{2,4,5} & \text{if $(i,j)=(2,4)$,}\\
(\la_5\Xi_{2,4}+\la_{2,5}\Xi_{2,5})/\la_{2,4,5} &\text{if 
$(i,j)=(2,5)$,}\\
0 &\text{if $(i,j)=(3,4)$,}\\
(\la_{5}\Xi_{2,3}-\la_{3}\Xi_{2,5})/\la_{1,3,5}&\text{if 
$(i,j)=(3,5)$,}\\
(-\la_{5}\Xi_{2,4}+\la_{4}\Xi_{2,5})/\la_{2,4,5}&\text{if 
$(i,j)=(4,5)$.}
\end{cases}
\]

The combinatorial type $\cT$ has a degeneration
of type $\cT'$ similar to $\cS'$, 
represented by the arrangement $\bar\A'$ 
in Figure \ref{fig:not selberg}. 
As in Example \ref{exm:Selberg}, 
the principal dependence is $(S,r)$, where $S=345$ and  $r=1$.
For $\cT$-nonresonant weights with $\la_S \neq 0$, 
the spectrum of $\GM_\bl(\cT',\cT)$ is contained in  
$\{0,\la_S\}$.
Calculations with the projection $\tau$ and the eigenspace decomposition 
of the endomorphism $\Omega_\bl(S,r)$ of $H^2(\cG)$ given in 
Example~\ref{exm:345} 
reveal that 
$\Omega_\bl(\cT',\cT)$ has eigenvalues 
$\la_S,\la_S,0$.
\end{exm}

\begin{exm}
The combinatorial type $\cS$ in Example \ref{exm:Selberg}
is a degeneration of the type $\cT$ in Example \ref{exm:not Selberg}.
The principal dependence of this degeneration is $(S,r)$, where $S=126$ and $r=2$. 
For $\cT$-nonresonant weights with $\la_S\neq 0$, the spectrum of $\Omega_\bl(\cS,\cT)$
is contained in $\{0,\la_S\}$. A calculation shows that the eigenvalues are $\la_S,0,0$.
It is interesting to note that $\la_S=\la_{1,2,6}=-\la_{3,4,5}$.
\end{exm}

\begin{ack}
This work was completed during the 
program ``Hyperplane Arrangements and Applications'' at the 
Mathematical Sciences Research Institute in Berkeley, California.  
We thank MSRI for its support and hospitality, and for providing a
stimulating mathematical environment.
\end{ack}

\bibliographystyle{amsalpha}

\end{document}